\documentclass[12pt]{article}

\usepackage{amsmath}
\usepackage{amsthm,amscd,amsfonts}
\usepackage{amssymb, upref, color}
\usepackage{amsmath,amssymb,amsthm}
\usepackage{color}
\usepackage[colorlinks]{hyperref}
\usepackage{graphicx}
\usepackage{epsf,epsfig,subfigure, verbatim}
\usepackage{latexsym,bm}

\numberwithin{equation}{section}

\theoremstyle{plain}

\long\def\@makefntext#1{\noindent #1}
\newskip\tabcentering \tabcentering=1000pt plus 1000pt minus 1000pt
\def\MCH#1#2{\setbox0=\hbox{\raise#1\hbox{#2}}\smash{\box0}}

\def\@evenfoot{}\def\@oddfoot{}

\def\@evenhead{\hbox to\textwidth{\footnotesize\rm\thepage \hfill
{\it  C. Zou \ \ Y. Xia \ M. Pinto\ J. Shi}}} 

\def\@oddhead{\hbox to \textwidth{\footnotesize{\it
Bounded Solutions and Topological Linearization  of
DEPCAGs} \hfill\thepage}}

\def\th#1{\vskip 1mm\noindent{\bf #1}\quad}

\def\proof{\vskip 1mm\noindent{\it Proof}\quad}

\catcode`@=12

\def\bc{\begin{center}}
\def\ec{\end{center}}
\def\no{\noindent}

\begin{document}
\date{}


\title{ Boundness and Linearisation of a class of differential equations with piecewise constant
argument 
}
\author{ \small Changwu Zou$^{1}$\footnote{ Changwu Zou was supported by the National Natural Science Foundation of China under Grant (No.11471027) and Foundation of Fujiang Province Education Department  under Grant (No. JAT160082). },
 \,\,  Yonghui Xia$^{2,3}$
 \footnote{ Corresponding author.  Yonghui Xia was supported by the National Natural
Science Foundation of China under Grant (No. 11671176 and No. 11271333), Natural
Science Foundation of Zhejiang Province under Grant (No. Y15A010022), Marie Curie Individual Fellowship within the European Community Framework Programme(MSCA-IF-2014-EF), the Scientific Research Funds of Huaqiao University and China Postdoctoral Science Foundation (No. 2014M562320). },
\,\,Manuel Pinto$^{4}$\footnote{Manuel Pinto was supported by FONDECYT Grant (No. 1120709 and No. 1170466 ).},\,\,Jinlin Shi$^{1}$,\,\,
\,\,Yuzhen Bai$^5$
\\
{\small 1.College of Mathematics and Computer Science, Fuzhou University, Fuzhou,
 350108 , China}\\
{\small \em zouchw@126.com (C. Zou)}\\
{\small 2. School of Mathematical Sciences, Huaqiao University, 362021, Quanzhou, Fujian, China.}\\
{\small\em xiadoc@163.com; xiadoc@hqu.edu.cn (Y.H.Xia)}\\
{\small 3. Department of Mathematics,
Zhejiang Normal University, Jinhua, 321004, China}\\
{\small 4. Departamento de Matematica, Universidad de Chile, Santiago, Chile}\\
{\small \em  pintoj.uchile@gmail.com(M. Pinto)    }\\
{\small 5. School of Mathematical Sciences, Qufu Normal University£¬
Qufu, 273165, P.R.China}\\
}
 \maketitle


\begin{center}
\begin{minipage}{140mm}
\begin{abstract}

The differential equations with piecewise constant argument (DEPCAs, for short) is a class of hybrid dynamical systems (combining continuous and discrete).  In this paper, under the assumption that the nonlinear term is partially unbounded, we study the bounded solution and global topological linearisation of a class of DEPCAs of general type. One of the purpose of this paper is to obtain a new criterion for the existence of a unique bounded solution, which improved the previous results. 
The other aim of this paper is to establish a generalized Grobman-Hartman result for the topological conjugacy between a nonlinear perturbation system and its linear system. The method is based on the new obtained criterion for bounded solution. The obtained results generalized and improved some previous papers. Some novel techniques are employed.

\end{abstract}



\end{minipage}
\end{center}
\section{\bf Introduction and Motivation}

 In this paper, we study the boundness and linearisation of a differential equations with piecewise constant argument of generalized type (for short, DEPCAGs). It takes the form
\begin{eqnarray}\label{sysnlz}
z'(t)=M(t)z(t)+M_0(t)z(\gamma(t))+h(t,z(t),z(\gamma(t))),
\end{eqnarray}
where $t\in \mathbb{R}, z(t)\in \mathbb{R}^{n}$, $M(t)$ and ${M_0}(t)$ are $n\times n$ matrices, $h:\mathbb{R}\times\mathbb{R}^{n}\times\mathbb{R}^{n}\rightarrow\mathbb{R}^n$ and $\gamma(t): \mathbb{R}\rightarrow\mathbb{R}$.

Usually, the authors studied the bounded solutions of perturbation nonlinear system under assumption that the perturbation nonlinear term is {\em bounded}.
 When the perturbation term is {\em unbounded}, it is difficult to study. In this paper, under the assumption that the nonlinear term $h(t,z(t),z(\gamma(t)))$ is partially unbounded, we study the bounded solution 
  of (\ref{sysnlz}). One of the purpose of this paper is to obtain a new criterion for the existence of a unique bounded solution, which improved the previous results (Theorem 5.3 in \cite{Coronel15}). Based on the new obtained criterion for bounded solution, we prove a generalized Grobamn-Hartman theorem to guarantee the conjugacy between the nonlinear system (\ref{sysnlz}) and its linear system.  This is another main purpose of this paper. The obtained results generalized and improved some previous papers. In fact, some novel techniques are employed in the proof.

Throughout this paper, we assume that the condition \textbf{(A)} holds:

There exist two constant sequences $\{t_i\}_{i\in\mathbb{Z}}$ and $\{\zeta_i\}_{i\in\mathbb{Z}}$ such that \\
\textbf{(A1)} $t_i<t_{i+1}$ and $t_i\leq\zeta_i\leq t_{i+1}$, $\forall i\in\mathbb{Z}$,\\
\textbf{(A2)} $t_i\rightarrow\pm\infty$ as $i\rightarrow\pm\infty,$\\
\textbf{(A3)} $\gamma(t)=\zeta_i$ for $t\in[t_i,t_{i+1})$,\\
\textbf{(A4)} there exists a constant $\theta>0$ such that $t_{i+1}-t_i\leq \theta, \forall i\in\mathbb{Z}.$

In particular, when $\gamma(t)=[t]$ or $\gamma(t)=2[\frac{(t+1)}{2}]$, system \eqref{sysnlz} is called the differential equations with piecewise constant argument (DEPCAs).

For DEPCAs  and DEPCAGs, many scholars study the continuity, boundedness, stability, existence of periodic or almost periodic solutions. One can refer to \cite{Aftabizadeh87,Akhmet08,Dai08,CastilloPinto15,ChiuPinto14,Pinto-Robledo15,VelozPinto15,Yuan02,Yuan97}. In particular, the bounded solutions of DEPCAs  and DEPCAGs were obtained in \cite{Akhmet07,ChiuPinto10,Coronel15,PintoRobledo15,PapaschinopoulosJMAA96}.
Among these works, Akhmet \cite{Akhmet07} obtained a set of sufficient conditions to guarantee the existence of a unique bounded solution by assuming that the linear system $z'(t)=M(t)z(t)$ in system \eqref{sysnlz} has an exponential dichotomy.  But if $M(t)=0$, then $z'(t)=M(t)z(t)$ can not admit an exponential dichotomy. It is possible that $z'(t)=M_0(t)z(\gamma(t))$ admits an exponential dichotomy, even if $M(t)=0$. In this case, the result in \cite{Akhmet07} is invalid. Later, Akhmet \cite{Akhmet12, Akhmet14} introduced the condition that the linear system with piecewise constant argument
\begin{eqnarray}\label{syslz}
z'(t)=M(t)z(t)+M_0(t)z(\gamma(t))
\end{eqnarray}
admits an exponential dichotomy. Under the assumption that linear system (\ref{syslz}) admits an exponential dichotomy and the nonlinear term $h(t,z(t),z(\gamma(t)))$ is bounded, Coronel et al \cite{Coronel15} proved that there exists a unique bounded solution of system (\ref{sysnlz}) (see \cite[Th. 5.3]{Coronel15}).
What happens if the nonlinear perturbed term $h(t,z(t),z(\gamma(t)))$ is unbounded? Does there exist a unique bounded solution?
This paper is devoted to answering this question. We prove that even if $h(t,z(t),z(\gamma(t)))$ is unbounded, system \eqref{sysnlz} has a unique bounded solution under some suitable conditions. We briefly summarize our result on bounded solution as follows:

\noindent {\bf Result 1}\, {\em Assume that system \eqref{syslz} admits an  exponential dichotomy and the nonlinear term $h(t,z(t),z(\gamma(t)))$ is Lipschitzian. If we further assume that there exist constants $r>0$  and $\mu>0$, such that
$$|h(t,z(t),z(\gamma(t)))|\leq r(|z(t)|+ |z(\gamma(t))|)+\mu,$$
then system \eqref{sysnlz} has a unique bounded solution under some conditions}.

\noindent {\bf Remark 1}\, We point out that $h(t,z(t),z(\gamma(t)))$ can be a polynomial of order one in $z(t)$ and $z(\gamma(t))$, which can be unbounded. For example, taking $h(t,z(t),z(\gamma(t)))= z(t)\sin t + z(\gamma(t))\cos t$. Thus our result improves Theorem 5.3 in \cite{Coronel15}. Certainly, we also generalize the results of Akhmet \cite{Akhmet07}.

Another purpose in this paper is to apply Result 1 to study the linearization of system \eqref{sysnlz} when the nonlinear term $h(t,z(t),z(\gamma(t)))$ is unbounded. Topological linearization is one of the most important research topics in the ordinary differential equations. A brief survey on topological linearization is presented as follows:


Since the Grobman-Hartman theorem was established by Hartman and Grobman \cite{Grobman62,Hartman63} in 1960's, many mathematician made contribution to this topic and made great progress in this theme. Most of the works focused on the autonomous systems. On the other hand, some mathematicians emphasized on the non-autonomous systems. Palmer \cite{Palmer73} proposed a version of Grobman-Hartman theorem for the non-autonomous ordinary differential equations in 1973.
For the ordinary differential equations, Barreira and Valls \cite{
Barreira06,Barreira11}, Jiang \cite{Jiang06,Jiang07}, 
Shi and Xiong \cite{Shi95} extended Palmer's result in various directions. For example, Shi and Xiong \cite{Shi95} reduced the conditions by assuming that the linear system partially admits an exponential dichotomy. Jiang \cite{Jiang06,Jiang07} reduced the condition by assuming that the linear system admits a generalized dichotomy. Barreira and Valls \cite{Barreira06,Barreira11} reduced the condition by assuming that the linear system admits a nonuniform exponential dichotomy. In addition, topological linearization of difference equations, functional differential equations and scalar reaction diffusion equations have been extensively studied. For examples, see  \cite{Castaneda15,Farkas,Kurzweil91,Lopez99, Lu,Papaschinopoulos94,Potzche08}. Another important work is the smooth linearization of $C^1$ hyperbolic mappings. One can refer to \cite{Bel2,ElBialy,Pugh,Sell,Sternberg1,Sternberg2,RS-JDE,RS-JDDE,ZWN-JDE,ZWN-JFA,ZWN-MA}.

In this paper, we focus on the topological linearization ($C^0$ linearization) of non-autonomous systems. As mentioned above, the topological linearization was extensively studied. However, seldom did the authors study the linearization problem of DEPCAGs. In 1996, Papaschinopoulos \cite{{Papaschinopoulos96}} generalized the topological linearization theorem to DEPCAs. And nineteen years later, Pinto and Robledo \cite{PintoRobledo15}  generalized the work of Papaschinopoulos to DEPCAGs. Under suitable conditions, they proved that the above nonlinear system \eqref{sysnlz} is topologically conjugated to its linear system \eqref{syslz}.
   They studied the linearization problem based on that the nonlinear terms in the systems are bounded. More specifically, the results in \cite{Papaschinopoulos96} and \cite{PintoRobledo15} require that $h(t,z(t),z(\gamma(t)))$ is bounded, ie., there exists a constant $\mu>0$ such that
   \[
   |h(t,z(t),z(\gamma(t)))|\leqslant \mu.
   \]
     However, in general, $h(t,z(t),z(\gamma(t)))$ can be unbounded. For example, taking $h(t,z(t),z(\gamma(t)))= z(t)\sin t + z(\gamma(t))\cos t$.  In this case, the results in \cite{Papaschinopoulos96} and \cite{PintoRobledo15} are not valid. In this paper, we prove that if $h(t,z(t),z(\gamma(t)))$ is unbounded, system \eqref{sysnlz} can also be topologically conjugated to system \eqref{syslz} as long as it has a proper structure. More precisely, we consider the following system with some proper structure
\begin{equation}\label{sysnlxy-phi}
\left\{\
\begin{array}{c}
x'(t)=A(t)x(t)+A_0(t)x(\gamma(t))+f(t,x(t),x(\gamma(t)))+\phi(t,y(t),y(\gamma(t))),\\
y'(t)=B(t)y(t)+B_0(t)y(\gamma(t))+g(t,x(t),x(\gamma(t)))+\psi(t,y(t),y(\gamma(t))),
\end{array}
\right.
\end{equation}
where $t\in \mathbb{R}, x(t)\in \mathbb{R}^{n_1}, y(t)\in \mathbb{R}^{n_2}$, $n_1+n_2=n$, $A(t)$, $A_0(t)$ are $n_1\times n_1$ matrices, $B(t)$, $B_0(t)$ are $n_2\times n_2$ matrices, $f:\mathbb{R}\times\mathbb{R}^{n_1}\times\mathbb{R}^{n_1}\rightarrow\mathbb{R}^{n_1}$,
$g:\mathbb{R}\times\mathbb{R}^{n_1}\times\mathbb{R}^{n_1}\rightarrow\mathbb{R}^{n_2}$, $\phi:\mathbb{R}\times\mathbb{R}^{n_1}\times\mathbb{R}^{n_1}\rightarrow\mathbb{R}^{n_1}$, and $\psi:\mathbb{R}\times\mathbb{R}^{n_2}\times\mathbb{R}^{n_2}\rightarrow\mathbb{R}^{n_2}$.

 In this paper, we assume that the nonlinear term is partially unbounded.
  We briefly summarize our second result on the global topological linearization of a class of DEPCAGs (\ref{sysnlxy-phi}) as follows:

\noindent {\bf Result 2}\, {\em
 Suppose that the linear system
\begin{equation}\label{syslxy}
\left\{\
\begin{array}{c}
x'(t)=A(t)x(t)+A_0(t)x(\gamma(t)),\\
y'(t)=B(t)y(t)+B_0(t)y(\gamma(t)),
\end{array}
\right.
\end{equation}
admits an exponential dichotomy. Assume that the nonlinear terms $f(t,x(t),x(\gamma(t))),\\ \,g(t,x(t),x(\gamma(t)))$ are Lipschitzian. If we further assume that there exist constants $\lambda$ and $\delta>0$ such that
$$|f(t,x(t),x(\gamma(t)))|\leq\lambda(|x(t)|+|x(\gamma(t))|),
\quad
 |g(t,x(t),x(\gamma(t)))|\leq\lambda(|x(t)|+|x(\gamma(t))|),$$
$$|\phi(t,y(t),y(\gamma(t)))|\leq \delta, \quad  |\psi(t,y(t),y(\gamma(t)))|\leq \delta. $$
Then system \eqref{sysnlxy-phi} is topologically conjugated to  system \eqref{syslxy} under proper conditions.
}

\noindent {\bf Remark 2}\, 
As you will see, the nonlinear terms $f(t,x(t),x(\gamma(t)))$ and  $g(t,x(t),x(\gamma(t)))$ can be possibly unbounded. For example, $f(t,x(t),x(\gamma(t)))$ and  $g(t,x(t),x(\gamma(t)))$  can be polynomials of order one in $x(t)$. In this case, the nonlinear term  of system \eqref{sysnlxy-phi} is unbounded in $x(t)$. For example, taking $f(t,x(t),x(\gamma(t)))=x(t)\sin t+x(\gamma(t))\cos t$. We see that the topological linearization can be realized. However, the results in \cite{{Papaschinopoulos96}} and \cite{PintoRobledo15} can not be applied to this case. In this sense, we extended the results in \cite{{Papaschinopoulos96}} and \cite{PintoRobledo15}.

It should be noted that if $\theta=0$, $A_0(t)=0$ and $B_0(t)=0$, \eqref{sysnlxy-phi} reduces to the ODE as follows.
\begin{equation}\label{sysnlxy-cor}
\left\{\
\begin{array}{c}
x'(t)=A(t)x(t)+f(t,x(t))+\phi(t,y(t)),\\
y'(t)=B(t)y(t)+g(t,x(t))+\psi(t,y(t)).
\end{array}
\right.
\end{equation}
Its linear system is
 \begin{equation}\label{syslxy-cor}
\left\{\
\begin{array}{c}
x'(t)=A(t)x(t),\\
y'(t)=B(t)y(t).
\end{array}
\right.
\end{equation}
Notice that if $\theta=0$, we can prove that if $|f(t,x(t)))|\leq\lambda(|x(t)|),
\quad
 |g(t,x(t))|\leq\lambda(|x(t)|)$, $|\phi(t,y(t))|\leq \delta, \quad  |\psi(t,y(t))|\leq \delta$,
then system \eqref{sysnlxy-cor} is topologically conjugated to system  \eqref{syslxy-cor}.

\th{Remark 3} Palmer \cite{Palmer73} proved that \eqref{ODE} is topologically conjugated to its linear part under the assumption that the nonlinear term is bounded. As you will see, the nonlinear terms $f$ and $g$ in system \eqref{sysnlxy-cor} can be unbounded. In this sense, we generalize and improve the main results in Palmer \cite{Palmer73}.

\noindent {\bf Remark 4}\, Some novel techniques are employed to prove our main result. Due to the unbounedness of the nonlinear terms, it is difficult to directly prove that the nonlinear system (\ref{sysnlxy-phi}) is topologically conjugated to the linear system (\ref{syslxy}). To overcome such difficulty, we should introduce the auxiliary system as follows
\begin{equation}\label{sysnlxy-fg}
\left\{\
\begin{array}{c}
x'(t)=A(t)x(t)+A_0(t)x(\gamma(t))+f(t,x(t),x(\gamma(t))),\\
y'(t)=B(t)y(t)+B_0(t)y(\gamma(t))+g(t,x(t),x(\gamma(t))).
\end{array}
\right.
\end{equation}
We first prove that system \eqref{syslxy} is topologically conjugated to system \eqref{sysnlxy-fg}. Secondly, we prove that system \eqref{sysnlxy-fg} and  system \eqref{sysnlxy-phi} are topologically conjugated. Then, by transition of topological conjugacy, system \eqref{syslxy} and  system  \eqref{sysnlxy-phi} are topologically conjugated.

The rest of this paper is organized as follows: In Section 2, we give some definitions, notation and preliminary lemmas. Our main results, Theorem 1 and Theorem 2, are stated in Section 3. The proof of Theorem 1 is  given in Section 4. The proof of Theorem 2 is very long and we divide the proofs into several Sections (see Sections  5-8).

\section{Preliminaries}

\no
\subsection{General assumptions}
\no

We introduce two groups of assumptions for Theorem 1 and Theorem 2, respectively. The first group is conditions \textbf{(B, C)} and the second group is conditions \textbf{($\mathfrak{B}, \mathfrak{C}$)}. In this paper, $|\cdot|$ denotes a vector norm or matrix  norm.

We assume that system \eqref{sysnlz} and \eqref{syslz} satisfy the following conditions.

Condition \textbf{(B)}:

\textbf{(B1)} The functions $M(t)$, $M_0(t)$ and $h(t,z(t),z(\gamma(t)))$ are locally integrable in $\mathbb{R}$.\\

\textbf{(B2)} There exists constants $r>0$, $\mu>0$ and $\ell>0$ such that for any $t\in \mathbb{R}$, $(t,z(t),z(\gamma(t)))$ and  $(t,\hat{z}(t),\hat{z}(\gamma(t)))$ $\in\mathbb{R}\times\mathbb{R}^{n}\times\mathbb{R}^{n}$,
$$
|h(t,z(t),z(\gamma(t)))|\leq r(|z(t)|+|z(\gamma(t))|)+\mu,
$$
and
$$
|h(t,z(t),z(\gamma(t)))-h(t,\hat{z}(t),\hat{z}(\gamma(t)))|\leq \ell\Big{(}|z(t)-\hat{z}(t)|+|z(\gamma(t))-\hat{z}(\gamma(t))|\Big{)}.
$$

We remark that if we further assume that $\ell\leq r$, $
|h(t,0,0)|\leq \mu
$, the Lipschitz condition in (B2) implies the first estimation $|h(t,z,y)|\leq r(|z|+|y|)+\mu$ in (B2).

Moreover, we introduce the following notation and condition \textbf{(C)}.
\begin{description}
\item[(i)] We define $I_i=[t_i,t_{i+1})$ for any $i\in\mathbb{Z}.$
\item[(ii)] For any $i\in\mathbb{Z}$ and $k\times k$ matrix $Q(t),$ we define
  $$\rho_i^+(Q)=\exp(\int_{t_i}^{\zeta_i}|Q(s)|ds) \quad \text{and} \quad  \rho_i^-(Q)=\exp(\int_{\zeta_i}^{t_{i+1}}|Q(s)|ds).$$
\end{description}

Condition \textbf{(C)}: There exists $0<\nu^+<1$ and $0<\nu^-<1$ such that the matrices $M(t)$ and $M_0(t)$ satisfy  properties:
$$
\sup\limits_{i\in\mathbb{Z}}\rho_i^+(M)\ln\rho_i^+(M_0)\leq\nu^+, \quad \sup\limits_{i\in\mathbb{Z}}\rho_i^-(M)\ln\rho_i^-(M_0)\leq\nu^-,
$$
and
\begin{equation}\label{defrho}
\begin{array}{c}
1\leq\rho(M)\triangleq\sup\limits_{i\in\mathbb{Z}}\rho_i^+(M)\rho_i^-(M)<+\infty.
\end{array}
\end{equation}

Therefore,
\begin{equation}\label{defalpha_0}
\begin{array}{c}
\rho_0(M)\triangleq\rho(M)^2(\frac{1+\nu^-}{1-\nu^+})>1.
\end{array}
\end{equation}

Now, we introduce conditions \textbf{($\mathfrak{ B}, \mathfrak{C}$)} for systems \eqref{sysnlxy-phi} and \eqref{syslxy}.

Condition \textbf{($\mathfrak{B}$)}:

\textbf{($\mathfrak{B}$1)} There exist constants $\beta>0$ and $\beta_0>0$ such that
$$\sup\limits_{t\in\mathbb{R}}|A(t)|\leq \beta, \quad \sup\limits_{t\in\mathbb{R}}|B(t)|\leq \beta,$$
$$\sup\limits_{t\in\mathbb{R}}|A_0(t)|\leq \beta_0, \quad \sup\limits_{t\in\mathbb{R}}|B_0(t)|\leq \beta_0.$$

\textbf{($\mathfrak{B}$2)} There exist constants $\delta>0$ and $\lambda>0$ such that for any $(t,x(t),x(\gamma(t)))\in\mathbb{R}\times\mathbb{R}^{n_1}\times\mathbb{R}^{n_1}$ and $(t,y(t),y(\gamma(t)))\in\mathbb{R}\times\mathbb{R}^{n_2}\times\mathbb{R}^{n_2},$
$$|f(t,x(t),x(\gamma(t)))|\leq \lambda(|x(t)|+|x(\gamma(t))|),$$
$$|g(t,x(t),x(\gamma(t)))|\leq \lambda(|x(t)|+|x(\gamma(t))|),$$
$$|\phi(t,y(t),y(\gamma(t)))|\leq\delta,$$
$$|\psi(t,y(t),y(\gamma(t)))|\leq\delta.$$

\textbf{($\mathfrak{B}$3)} There exists constant $\omega>0$ such that for any $(t,x_1(t),x_1(\gamma(t)))$, $(t,x_2(t),x_2(\gamma(t)))$ $\in\mathbb{R}\times\mathbb{R}^{n_1}\times\mathbb{R}^{n_1}$ and  $(t,y_1(t),y_1(\gamma(t)))$, $(t,y_2(t),y_2(\gamma(t)))$ $\in\mathbb{R}\times\mathbb{R}^{n_2}\times\mathbb{R}^{n_2}$,
$$
\begin{array}{ccc}
  & &|f(t,x_{1}(t),x_{1}(\gamma(t)))-f(t,x_{2}(t),x_{2}(\gamma(t)))|\\
  &\leq& \omega\Big{(}|x_{1}(t)-x_{2}(t)|+|x_{1}(\gamma(t))-x_{2}(\gamma(t))|\Big{)},
\end{array}
$$
$$
\begin{array}{ccc}
  & &|g(t,x_{1}(t),x_{1}(\gamma(t)))-g(t,x_{2}(t),x_{2}(\gamma(t)))|\\
  &\leq& \omega\Big{(}|x_{1}(t)-x_{2}(t)|+|x_{1}(\gamma(t))-x_{2}(\gamma(t))|\Big{)},
\end{array}
$$
$$
\begin{array}{ccc}
  & &|\phi(t,y_{1}(t),y_{1}(\gamma(t)))-\phi(t,y_{2}(t),y_{2}(\gamma(t)))|\\
  &\leq& \omega\Big{(}|y_{1}(t)-y_{2}(t)|+|y_{1}(\gamma(t))-y_{2}(\gamma(t))|\Big{)},
\end{array}
$$
and
$$
\begin{array}{ccc}
  & &|\psi(t,y_{1}(t),y_{1}(\gamma(t)))-\psi(t,y_{2}(t),y_{2}(\gamma(t)))|\\
  &\leq& \omega\Big{(}|y_{1}(t)-y_{2}(t)|+|y_{1}(\gamma(t))-y_{2}(\gamma(t))|\Big{)}.
\end{array}
$$

Condition \textbf{($\mathfrak{C}$)}: There exist $0<\nu^+<1$ and $0<\nu^-<1$ such that matrices $A(t)$, $A_0(t)$, $B(t)$ and $B_0(t)$ satisfy following properties:
$$
\sup\limits_{i\in\mathbb{Z}}\rho_i^+(A)\ln\rho_i^+(A_0)\leq\nu^+, \quad \sup\limits_{i\in\mathbb{Z}}\rho_i^-(A)\ln\rho_i^-(A_0)\leq\nu^-,
$$
$$
\sup\limits_{i\in\mathbb{Z}}\rho_i^+(B)\ln\rho_i^+(B_0)\leq\nu^+, \quad \sup\limits_{i\in\mathbb{Z}}\rho_i^-(B)\ln\rho_i^-(B_0)\leq\nu^-.
$$
Note that \textbf{($\mathfrak{B}$1)} and \textbf{(A4)} imply that
\begin{equation}\label{defrho_AB}
\begin{array}{c}
1\leq\rho(A)\triangleq\sup\limits_{i\in\mathbb{Z}}\rho_i^+(A)\rho_i^-(A)<+\infty  \quad \textrm{and} \quad
1\leq\rho(B)\triangleq\sup\limits_{i\in\mathbb{Z}}\rho_i^+(B)\rho_i^-(B)<+\infty.
\end{array}
\end{equation}
Thus,
\begin{equation}\label{defalpha_0AB}
\begin{array}{c}
\rho_0(A)\triangleq\rho^2(A)(\frac{1+\nu^-}{1-\nu^+})>1  \quad \textrm{and}\quad
\rho_0(B)\triangleq\rho^2(B)(\frac{1+\nu^-}{1-\nu^+})>1.
\end{array}
\end{equation}

Throughout the rest of the paper, we assume that conditions \textbf{(A, B, C, $\mathfrak{B}, \mathfrak{C}$)} hold.

\subsection{Notation of solutions for DEPCAGs}
\no

The notion of solutions for DEPCAGs was introduced in \cite{Aftabizadeh87,Akhmet14,ChiuPinto10,CookeWiener89,Coronel15,Wiener93}.
\th{Definition 2\ (Solutions of a DEPCAG)}
A continuous function $z(t)$ is a solution of system \eqref{sysnlz} or system \eqref{syslz} on $\mathbb{R}$ if:
\begin{description}
  \item[(i)] The derivative $z'(t)$ exists at each point $t\in\mathbb{R}$ with the possible exception of points $t_i,i\in\mathbb{Z}$, where the one side derivative exists;
  \item[(ii)] The equation is satisfied for $z(t)$ on each interval $(t_i,t_{i+1})$ and it holds for the right derivative of $z(t)$ at $t_i$.
\end{description}

\subsection{Transition matrices}
\no

In this subsection, we introduce some notation associated with solutions of a class of DEPCAGs.

Let $\Phi(t)$ be the fundamental matrix of system $x'=M(t)x$ with $\Phi(0)=I$.
For any $t\in I_j$, $\tau\in I_i$, $s\in \mathbb{R}$, we introduce the following notations adopting from \cite{Coronel15,PintoJDEQ11,PintoRobledo15}:
$$\Phi(t,s)=\Phi(t)\Phi^{-1}(s),$$
$$J(t,\tau)=I+\int_\tau^t\Phi(\tau,s)M_0(s)ds,$$
$$E(t,\tau)=\Phi(t,\tau)+\int_\tau^t\Phi(t,s)M_0(s)ds=\Phi(t,\tau)J(t,\tau).$$

We define backward and forward products of a set of $k\times k$ matrices $\mathcal{Q}_i (i=1,\ldots,m)$ as follows:
$$
\prod\limits_{i=1}^{\leftarrow m}\mathcal{Q}_i=\left\{ \begin{array}{cc}
\mathcal{Q}_m\cdots\mathcal{Q}_2\mathcal{Q}_1, &\quad\text{if} \quad m\geqslant  1,\\
I, &\quad \text{if} \quad  m<1,
\end{array}
\right.
$$
and
$$
\prod\limits_{i=1}^{\rightarrow m}\mathcal{Q}_i=\left\{ \begin{array}{cc}
\mathcal{Q}_1\mathcal{Q}_2\cdots\mathcal{Q}_m, &\quad\text{if} \quad m\geqslant  1,\\
I, &\quad \text{if} \quad  m<1.
\end{array}
\right.
$$

If $J(t,s)$ is nonsingular, we could define the transition matrix $Z(t,s)$ of system \eqref{syslz} as follows:\\
if $t>\tau$,
\begin{eqnarray*}
&&Z(t,\tau)\\
&=&E(t,\zeta_j)E(t_j,\zeta_j)^{-1}\prod\limits_{r=i+2}^{\leftarrow j}\Big{(}E(t_r,\gamma(t_{r-1}))E(t_{r-1},\gamma(t_{r-1}))^{-1}\Big{)}E(t_{i+1},\gamma(\tau))E(\tau,\gamma(\tau))^{-1},\nonumber
\end{eqnarray*}
if $t<\tau$,
\begin{eqnarray*}
&&Z(t,\tau)\\
&=&E(t,\zeta_j)E(t_{j+1},\zeta_j)^{-1}\prod\limits_{r=j+1}^{\rightarrow i-1}\Big{(}E(t_{r},\gamma(t_{r}))E(t_{r+1},\gamma(t_{r}))^{-1}\Big{)}E(t_{i},\gamma(\tau))E(\tau,\gamma(\tau))^{-1}.\nonumber
\end{eqnarray*}

Through simple calculations, we obtain $Z(t,\tau)Z(\tau,s)=Z(t,s)$ and $Z(t,s)=Z(s,t)^{-1}$. Since $E(\tau,\tau)=I$ and $\frac{\partial E(t,\tau)}{\partial t}=M(t)E(t,\tau)+M_0(t)$,
we have
\begin{equation*}
\frac{\partial Z(t,\tau)}{\partial t}=M(t)Z(t,\tau)+M_0(t)Z(\gamma(t),\tau).
\end{equation*}
Thus, $Z(t,\tau)$ is a solution of system \eqref{syslz}.

\subsection{Formulas of solutions for DEPCAGs}
\no

To introduce the formulas of solutions, we first state the following important lemma.
\th{Lemma 2.1(\cite{PintoJDEQ11}, Lemma 4.3)}
\emph{Assume that conditions \textbf{(A,B,C)} are fulfilled, then $J(t,s)$ is nonsingular for any $t,s\in \bar{I}_r$ and the matrices
$Z(t,s)$ and $Z(t,s)^{-1}$ are well defined for any $t, s\in\mathbb{R}$. If $t,s\in \bar{I}_r$, then
\begin{equation*}
|\Phi(t,s)|\leq\rho(M),
\end{equation*}
\begin{equation*}
|Z(t,s)|\leq\rho_0(M),
\end{equation*}
where $\rho(\cdot)$ is defined in \eqref{defrho} and $\rho_0(\cdot)$ is defined in \eqref{defalpha_0}.}

We remark that Lemma 2.1 ensures the continuity of solutions of system \eqref{sysnlz} on $\mathbb{R}$. We introduce the following formulas for DEPCAGs.

\th{Proposition 2.1 (\cite{PintoJDEQ11}, p.239)}\emph{For any $t\in I_j$, $\tau\in I_i,$ the solution of system \eqref{syslz} with $x(\tau)=\xi$ is defined on $\mathbb{R}$ and is given by
\begin{equation}\label{sollz}
z(t)=Z(t,\tau)\xi.
\end{equation}  }

\th{Proposition 2.2 (\cite{PintoJDEQ11}, Th 3.3)}\emph{For any $t\in I_j$, $\tau\in I_i$ and $t>\tau$, the solution of system \eqref{sysnlz} with $z(\tau)=\xi$ is defined on $\mathbb{R}$ and is given by
\begin{eqnarray}\label{solnlz}
z(t)&=&Z(t,\tau)\xi+\int_\tau^{\zeta_i}Z(t,\tau)\Phi(\tau,s)h(s)ds
+\sum\limits_{r=i+1}^j\int_{t_r}^{\zeta_r}Z(t,t_r)\Phi(t_r,s)h(s)ds \nonumber\\
& & +\sum\limits_{r=i}^{j-1}\int_{\zeta_r}^{t_{r+1}}Z(t,t_{r+1})\Phi(t_{r+1},s)h(s)ds
+\int_{\zeta_j}^{t}\Phi(t,s)h(s)ds,        
\end{eqnarray}
where $h(s)=h(s,z(s),z(\gamma(s)))$.  }

\th{Remark 2.1}If $t<\tau$, one could obtain the solution formula by replacing $\sum\limits_{r=i+1}^j$ and $\sum\limits_{r=i}^{j-1}$ by $\sum\limits_{r=j+1}^i$ and $\sum\limits_{r=j}^{i-1}$, respectively.

\subsection{Subsystems of System \eqref{sysnlxy-phi}}
\no

For convenience, consider the following subsystems of system \eqref{sysnlxy-phi}:
\begin{equation}\label{sysnlx}
x'(t)=A(t)x(t)+A_0(t)x(\gamma(t))+f(t,x(t),x(\gamma(t)))+\phi(t,y(t),y(\gamma(t))),
\end{equation}
\begin{equation}\label{sysnly}
 y'(t)=B(t)y(t)+B_0(t)y(\gamma(t))+g(t,x(t),x(\gamma(t)))+\psi(t,y(t),y(\gamma(t))),
\end{equation}
and subsystems of system \eqref{syslxy}:
\begin{equation}\label{syslx}
x'(t)=A(t)x(t)+A_{0}(t)x(\gamma(t)),
\end{equation}
\begin{equation}\label{sysly}
 y'(t)=B(t)y(t)+B_{0}(t)y(\gamma(t)).
\end{equation}

Let $\Phi_1(t)$ be the fundamental matrix of system $x'=A(t)x$ with $\Phi_1(0)=I$, and $\Phi_2(t)$ be the fundamental matrix of system $y'=B(t)y$ with $\Phi_2(0)=I$.

For any $t\in I_j$, $\tau\in I_i$, $s\in \mathbb{R}$, similar to $\Phi(t,s), J(t,\tau)$, and $E(t,\tau)$ in subsection 2.3, we could define
$$\Phi_k(t,s), \quad J_k(t,\tau)\quad  \text{and} \quad E_k(t,\tau),\quad k=1,2.$$

If $J_k(t,s)$ $(k=1,2)$ is nonsingular,  we could define the transition matrices  $Z_1(t,s)$ and $Z_2(t,s)$ of subsystems \eqref{syslx} and \eqref{sysly}, respectively. Moreover, we could verify that $Z_1(t,\tau)$ and $Z_2(t,\tau)$ are solutions of subsystems \eqref{syslx} and \eqref{sysly}, respectively.

\subsection{$\alpha$-exponential dichotomy and Green function}
\no

Now we introduce the definition of exponential dichotomy for a DEPCAG. In this paper, we adopt the following definition from Akhmet \cite{Akhmet12, Akhmet14}.
\th{Definition 3\ ($\alpha$-exponential dichotomy for a DEPCAG) }
The linear system \eqref{syslz} has an $\alpha$-exponential dichotomy on $\mathbb{R}$ if there exist a projection $P$, constants $K\geqslant  1$ and $\alpha>0$ such that the transition matrix $Z(t,s)$ of system \eqref{syslz} satisfies
\begin{equation*}
|Z_P(t,s)|\leq Ke^{-\alpha|t-s|},
\end{equation*}
where $Z_P(t,s)$ is defined by
\begin{equation*}
Z_P(t,s)=\left\{ \begin{array}{cc}
Z(t,0)PZ(0,s), & t\geqslant  s,\\
-Z(t,0)(I-P)Z(0,s), &  s>t.
\end{array}
\right.
\end{equation*}

For convenience, we define the Green function corresponding to system \eqref{sysnlz} which was introduced in \cite{PintoRobledo15,Coronel15}.
Given $t\in(\zeta_j,t_{j+1})$,
$$
\tilde{G}(t,s)= \left\{ \begin{array}{ccc}
Z_p(t,t_r)\Phi(t_r,s), & \textrm{ if } & s\in[t_r,\zeta_r) \textrm{ for any } r\in\mathbb{Z},\\
Z_p(t,t_{r+1})\Phi(t_{r+1},s) & \textrm{ if } & s\in[\zeta_r,t_{r+1}) \textrm{ for any } r\in\mathbb{Z}\setminus\{j\},\\
\Phi(t,s) & \textrm{ if } & s\in[\zeta_j,t),\\
 0 & \textrm{ if } & s\in[t,t_{j+1}],
\end{array}
\right.
$$
and if $t\in[t_j,\zeta_j],$
$$
\tilde{G}(t,s)= \left\{ \begin{array}{ccc}
Z_p(t,t_r)\Phi(t_r,s) & \textrm{ if } & s\in[t_r,\zeta_r) \textrm{ for any } r\in\mathbb{Z}\setminus\{j\},\\
Z_p(t,t_{r+1})\Phi(t_{r+1},s) & \textrm{ if } & s\in[\zeta_r,t_{r+1}) \textrm{ for any } r\in\mathbb{Z},\\
0 & \textrm{ if } & s\in[t_j,t),\\
-\Phi(t,s) & \textrm{ if } & s\in[t,\zeta_j),
\end{array}
\right.
$$

We denote $\tilde{G}_1(t,s)=\tilde{G}(t,s)$ for $t\geqslant s$ and  $\tilde{G}_2(t,s)=-\tilde{G}(t,s)$ for $t<s$.

\subsection{Condition \textbf{($\mathfrak{D}$)}}
\no

For convenience, we apply the following condition in our second result to replace the condition that system \eqref{syslxy} has an $\alpha$-exponential dichotomy.

Condition \textbf{($\mathfrak{D}$)}: There exist constants $K\geqslant  1$ and $\alpha>0$ such that
$$
|Z_1(t,s)|\leq e^{-\alpha(t-s)}, \quad t\geqslant  s \quad \textrm{and} \quad  |Z_2(t,s)|\leq Ke^{\alpha(t-s)}, \quad s>t.
$$

It is clear that condition \textbf{($\mathfrak{D}$)} is equivalent to assume that system \eqref{syslxy} has an $\alpha$-exponential dichotomy by taking $K=1$ in the first inequality. We point out that this assumption is natural. 
In fact, we can get the  inequality by taking  another equivalent norm  or supposing the following conditions:
$$
\frac{d|x(t)|'}{dt}|_{\eqref{syslxy}}\leq -2\alpha|x(t)|^2, \quad |f(t,x(t),x(\gamma(t))|\leq \frac{\alpha}{2}|x|.
$$

\subsection{Topological conjugacy}
\no

The notion of topological equivalence and topological conjugacy can be found in \cite{Palmer73,Palmer75,PintoRobledo15,XiaLi13}.
\th{Definition 1\ (Topological conjugacy)}
A continuous function $H:\mathbb{R}\times\mathbb{R}^n\rightarrow\mathbb{R}^n$ is topological equivalence between system \eqref{sysnlz} and \eqref{syslz} if following conditions hold:
\begin{description}
  \item[(i)] for each $t\in\mathbb{R}$, $H(t,z)$ is a homeomorphism of $\mathbb{R}^n$,
  \item[(ii)] $H(t,z) \rightarrow \infty$  as $z \rightarrow \infty$ uniformly with respect to t,
  \item[(iii)] if $z(t)$ is a solution of system \eqref{sysnlz}, then $H(t,z(t))$ is a solution of system \eqref{syslz}.
\end{description}
In addition, the function $L(t,z)=H^{-1}(t,z)$ has properties (i)-(iii) also.

If such a map $H$ exists, then system \eqref{sysnlz} and \eqref{syslz} are called topologically conjugated.

\subsection{Some lemmas}
\no

\th{Lemma \ 2.2 (\cite{PintoRobledo15}, Proposition 3)}
\emph{If  system \eqref{syslz} has an $\alpha$-exponential dichotomy on $\mathbb{R}$, then $\tilde{G}$ satisfies
$$ |\tilde{G}(t,s)|\leq K\rho^\ast(M) e^{-\alpha|t-s|},$$
where $\rho^\ast(M)=\rho(M)e^{\alpha\theta}$, $\rho(M)$ is defined in \eqref{defrho} and $\theta$ is in \textbf{(A4)}.  }

From Lemma 2.2, we have that
\begin{equation}\label{tilG12}
|\tilde{G}_1(t,s)|\leq K\rho^\ast(M) e^{-\alpha(t-s)} \quad \text{for} \quad t\geqslant s,\quad  |\tilde{G}_2(t,s)|\leq K\rho^\ast(M) e^{-\alpha(s-t)}\quad \text{for} \quad t< s.
\end{equation}

\th{Lemma \ 2.3 (\cite{PintoRobledo15}, Lemma 2.3)}
\emph{If  system \eqref{syslz}   has an $\alpha$-exponential dichotomy on $\mathbb{R}$,  then the unique solution bounded on $\mathbb{R}$ is the null solution.}

\section{Main Results}
\no
Now we are in a position to state our main results.

\th{Theorem 1} \emph{If conditions \textbf{(A,B,C)} hold and system \eqref{syslz} has an $\alpha$-exponential dichotomy with constant $K\geqslant1$ and $\alpha>0$, further assume that
\begin{equation}\label{eq3.1}
8Kl\rho^{\ast}(M)\alpha^{-1}\leq 1, \quad  4Kr\rho^{\ast}(M)\alpha^{-1}\leq 1,
\end{equation}
where $\rho^{\ast}(M)$ is defined in Lemma 2.2, then system \eqref{sysnlz} has a unique solution bounded on $\mathbb{R}$ which can be represented as follows
$$
z(t) = \int_{-\infty}^t\tilde{G}_1(t,s)h(s,z(s),z(\gamma(s)))ds - \int_t^{+\infty}\tilde{G}_2(t,s)h(s,z(s),z(\gamma(s)))ds,
$$
and
$$
|z(t)| \leq 2 K\mu\tilde{\rho}(M)(\alpha-4rK\tilde{\rho}(M))^{-1} \triangleq \sigma.
$$
}

We remark that if system \eqref{sysnlz} reduces to ODE, that is,
\[
z'(t)=M(t)z(t) +h(t,z(t)),
\]
Theorem 1 is valid for ODE.

\th{Theorem\ 2}\emph{If conditions \textbf{(A,$\mathfrak{B}$,$\mathfrak{C}$,$\mathfrak{D}$)} hold,  further assume that
\begin{equation}\label{eq3.2}
8K\tilde{\rho}(A)\omega\alpha^{-1} < 1, \quad 8K\tilde{\rho}(B)\omega\alpha^{-1} <1,
\end{equation}
\begin{equation}\label{eq3.2+1}
16K\tilde{\rho}(A)\lambda\alpha^{-1} < 1, \quad 16K\tilde{\rho}(B)\lambda\alpha^{-1} <1,
\end{equation}
\begin{equation}\label{eq3.2+2}
\alpha_0=\alpha-2\omega\tilde{\rho}(A)e^{\alpha\theta}>0,
\end{equation}
\begin{equation}\label{eq3.3}
F(\ell,\theta)(\beta_0+\ell)\theta = \upsilon < 1,
\end{equation}
where $F(\ell,\theta) = \frac{e^{(\beta+\ell)\theta}-1}{(\beta+\ell)\theta}$, $\tilde{\rho}(\cdot)=\max(\rho(\cdot)\rho_0(\cdot), \rho(\cdot)e^{\alpha\theta})$, $\rho(\cdot)$ is defined in \eqref{defrho_AB}, $\rho_0(\cdot)$ is defined in \eqref{defalpha_0AB},
 then system \eqref{sysnlxy-phi} is topologically conjugated to system \eqref{syslxy}.   }

If system \eqref{sysnlz} reduces to ODE, that is,
\begin{equation}\label{ODE}
z'(t)=M(t)z(t) +h(t,z(t)).
\end{equation}
And \eqref{sysnlxy-phi} reduces to the ODE system \eqref{sysnlxy-cor}. System \eqref{syslxy-cor} is its linear system.
Notice that for the ODE system, $\theta=0$, then $\rho^{\ast}(A)=1$, $\rho^{\ast}(B)=1$. Thus, Theorem 2 reduces to

\th{Corollary\ 1}  \emph{Assume that system \eqref{syslxy-cor} has an $\alpha$-exponential dichotomy with constant $K\geqslant1$ and $\alpha>0$. If $|f(t,x(t)))|\leq\lambda(|x(t)|),
\quad
 |g(t,x(t))|\leq\lambda(|x(t)|)$, $|\phi(t,y(t))|\leq \delta, \quad  |\psi(t,y(t))|\leq \delta$, and further assume that
$$
8K\omega\alpha^{-1} < 1, \quad 8K\lambda\alpha^{-1} <1,
$$
then system \eqref{sysnlxy-cor} is topologically conjugated to  system  \eqref{syslxy-cor}. }


\section{The Proof of Theorem 1 }
\no

To prove Theorem 1, we first introduce the following  lemma.

\th{Lemma 4.1}\emph{ If $t>\zeta_i$ and $z(t)$ is a bounded solution of system \eqref{sysnlz}, then
\begin{eqnarray*}
I&\triangleq&\sum_{r=-\infty}^{i}\int_{t_r}^{\zeta_r}Z(t,0)P Z(0,t_r)\Phi(t_r,s)h(s,z(s),z(\gamma(s)))ds
\end{eqnarray*}
is convergent. }
\proof From $t_r\leq\zeta_r$, $t>\zeta_i$,  \textbf{(B)} and \eqref{tilG12},   we have
\begin{eqnarray*}
|I|&\leq&\int_{-\infty}^{t}|Z(t,0)P Z(0,t_r)\Phi(t_r,s)h(s,z(s),z(\gamma(s)))|ds \\
   &\leq&\int_{-\infty}^{t}|\tilde{G}_1(t,s)|{r(|z(s)|+|z(\gamma(s))|)+\mu}ds   \\
  &\leq& \int_{-\infty}^{t}K\rho^\ast(M) e^{-\alpha(t-s)} (2r |z|+\mu)ds  \\
  &=& K\rho^\ast(M)\alpha^{-1}(2r |z|+\mu).
\end{eqnarray*}
Since $z(s)$ is a bounded solution, $I$  is convergent.    $\Box$

\th{The proof of Theorem 1:}
\proof For $\sigma$ defined in Theorem 1, denote
$$\Omega=\{\varphi(t)|\varphi:\mathbb{R}\rightarrow \mathbb{R}^n \text{ is continuous and } |\varphi(t)|\leq\sigma\},$$
and
$$
W=\{\varphi(t)|\varphi:\mathbb{R}\rightarrow\mathbb{R}^n\text{ is continuous and }\|\varphi\| < \infty\}.$$
It is easy to see that $W$ is a Banach space and $\Omega$ is a closed subset of $W$.

Suppose that  $t\in{[\zeta_j,t_{j+1})}, 0\in{[t_i,\zeta_i]},j>i.$
For any $\varphi(t)\in\Omega$, define the map $T:\Omega\rightarrow W$ as follows
\begin{align}
T\varphi(t) \
&= \sum_{r=-\infty}^{j}\int_{t_r}^{\zeta_r}Z(t,0)PZ(0,t_r)\Phi(t_r,s)h(s,\varphi(s),\varphi(\gamma(s)))ds \notag\\
&+ \sum_{r=-\infty}^{j-1}\int_{\zeta_r}^{t_{r+1}}Z(t,0)PZ(0,t_{r+1})\Phi(t_{r+1},s)h(s,\varphi(s),\varphi(\gamma(s)))ds \notag\\
&- \sum_{r=j+1}^{+\infty}\int_{t_r}^{\zeta_r}Z(t,0)(I-P)Z(0,t_r)\Phi(t_r,s)h(s,\varphi(s),\varphi(\gamma(s)))ds \notag\\
&- \sum_{r=j}^{+\infty}\int_{\zeta_r}^{t_{r+1}}Z(t,0)(I-P)Z(0,t_{r+1})\Phi(t_{r+1},s)h(s,\varphi(s),\varphi(\gamma(s)))ds \notag\\
&+ \int_{\zeta_j}^{t}\Phi(t,s)h(s,\varphi(s),\varphi(\gamma(s)))ds \notag \\
&=\int_{-\infty}^t\tilde{G}_1(t,s)h(s,\varphi(s),\varphi(\gamma(s)))ds - \int_t^{+\infty}\tilde{G}_2(t,s)h(s,\varphi(s),\varphi(\gamma(s)))ds. \notag
\end{align}

To prove the existence and uniqueness of bounded solution, we divide it into two steps.

\noindent {\bf Step 1}\, We prove that the map $T$ has a unique fixed point by contraction principle.

Due to \eqref{tilG12} and \textbf{(B2)}, we get
\begin{align}
|T\varphi(t)| \
&\leq \int_{-\infty}^{t}Ke^{-\alpha(t-s)}\tilde{\rho}(M)(r|z(t)|+r|z(\gamma(t))|+\mu)ds  \notag\\
&\quad + \int_{t}^{+\infty}Ke^{\alpha(t-s)}\tilde{\rho}(M)(r|z(t)|+r|z(\gamma(t))|+\mu)ds \notag\\
&\leq [K\tilde{\rho}(M)(\mu+2r\sigma)+K\tilde{\rho}(M)(\mu+2r\sigma)]\alpha^{-1} \notag\\
&= 2K\tilde{\rho}(M)\alpha^{-1}(\mu+2r\sigma) \notag\\
&=\sigma. \notag
\end{align}

Therefore $T\varphi\in\Omega$ and $T$ is a map from $\Omega$ to $\Omega$.

For any $\varphi_1(t),\varphi_2(t)\in\Omega$, from \eqref{tilG12} and \textbf{(B2)}
we have
\begin{align}
|T\varphi_1(t)-T\varphi_2(t)| \
&=\Big{|}\int_{-\infty}^{t}\tilde{G}_1(t,s)[h(s,\varphi_1(s),\varphi_1(\gamma(s)))-h(s,\varphi_2(s),\varphi_2(\gamma(s)))]ds\notag\\
&\quad+\int_{t}^{+\infty}\tilde{G}_2(t,s)[h(s,\varphi_1(s),\varphi_1(\gamma(s)))-h(s,\varphi_2(s),\varphi_2(\gamma(s)))]ds \Big{|}\notag\\
&\leq \int_{-\infty}^{t}K\tilde{\rho}(M)e^{-\alpha(t-s)}l(|\varphi_1(s)-\varphi_2(s)|+|\varphi_1(\gamma(s))-\varphi_2(\gamma(s))|)\notag\\
&\quad+\int_{t}^{+\infty}K\tilde{\rho}(M)e^{\alpha(t-s)}l(|\varphi_1(s)-\varphi_2(s)|+|\varphi_1(\gamma(s))-\varphi_2(\gamma(s))|)ds\notag\\
&\leq 2Kl\tilde{\rho}(M)\alpha^{-1}\|\varphi_1-\varphi_2\| + 2Kl\tilde{\rho}(M)\alpha^{-1}\|\varphi_1-\varphi_2\| \notag\\
&\leq \frac{1}{2}\|\varphi_1-\varphi_2\|. \notag
\end{align}

Thus $T$ is a contraction map in $\Omega$. By the contraction map principle, there exists a unique $\varphi_0(t) \in \Omega$ such that
\begin{align}
\varphi_0(t)&= T\varphi_0(t) \notag\\
&= \sum_{r=-\infty}^{j}\int_{t_r}^{\zeta_r}Z(t,0)PZ(0,t_r)\Phi(t_r,s)h(s,\varphi_0(s),\varphi_0(\gamma(s)))ds \notag\\
&+ \sum_{r=-\infty}^{j-1}\int_{\zeta_r}^{t_{r+1}}Z(t,0)PZ(0,t_{r+1})\Phi(t_{r+1},s)h(s,\varphi_0(s),\varphi_0(\gamma(s)))ds \notag\\
&- \sum_{r=j+1}^{+\infty}\int_{t_r}^{\zeta_r}Z(t,0)(I-P)Z(0,t_r)\Phi(t_r,s)h(s,\varphi_0(s),\varphi_0(\gamma(s)))ds \notag\\
&- \sum_{r=j}^{+\infty}\int_{\zeta_r}^{t_{r+1}}Z(t,0)(I-P)Z(0,t_{r+1})\Phi(t_{r+1},s)h(s,\varphi_0(s),\varphi_0(\gamma(s)))ds \notag\\
&+ \int_{\zeta_j}^{t}\Phi(t,s)h(s,\varphi_0(s),\varphi_0(\gamma(s)))ds. \notag
\end{align}
Furthermore, it is easy to check that $\varphi_0(t)$ is a solution of system \eqref{sysnlz}.

\noindent {\bf Step 2}\, We prove the uniqueness of the bounded solution.
That is, we prove that $\varphi_0(t)$ is the unique bounded solution of system \eqref{sysnlz}.
In fact, suppose that $\varphi_1(t)$ is another bounded solution of system \eqref{sysnlz}, by Proposition 2.2, we get
\begin{align}
\varphi_1(t) \
&=Z(t,0)\varphi_1(0)+\int_{0}^{\zeta_i}Z(t,0)\Phi(0,s)h(s,\varphi_1(s),\varphi_1(\gamma(s)))ds\notag\\
&\quad+\sum_{r=i+1}^{j}\int_{t_r}^{\zeta_r}Z(t,t_r)\Phi(t_r,s)h(s,\varphi_1(s),\varphi_1(\gamma(s)))ds\notag\\
&\quad+\sum_{r=i}^{j-1}\int_{\zeta_r}^{t_{r+1}}Z(t,t_{r+1})\Phi(t_{r+1},s)h(s,\varphi_1(s),\varphi_1(\gamma(s)))ds\notag\\
&\quad+\int_{\zeta_j}^{t}\Phi(t,s)h(s,\varphi_1(s),\varphi_1(\gamma(s)))ds\notag\\
&=Z(t,0)\{\varphi_1(0)+\int_{0}^{\zeta_i}\Phi(0,s)h(s,\varphi_1(s),\varphi_1(\gamma(s)))ds\notag\\
&\quad+\sum_{r=i+1}^{j}P\int_{t_r}^{\zeta_r}Z(0,t_r)\Phi(t_r,s)h(s,\varphi_1(s),\varphi_1(\gamma(s)))ds\notag\\
&\quad+\sum_{r=i}^{j-1}P\int_{\zeta_r}^{t_{r+1}}Z(0,t_{r+1})\Phi(t_{r+1},s)h(s,\varphi_1(s),\varphi_1(\gamma(s)))ds\}\notag\\
&\quad+\int_{\zeta_j}^{t}\Phi(t,s)h(s,\varphi_1(s),\varphi_1(\gamma(s)))ds \notag\\
&\quad+Z(t,0)\{\sum_{r=i+1}^{j}(I-P)\int_{t_r}^{\zeta_r}Z(0,t_r)\Phi(t_r,s)h(s,\varphi_1(s),\varphi_1(\gamma(s)))ds\notag\\
&\quad+\sum_{r=i}^{j-1}(I-P)\int_{\zeta_r}^{t_{r+1}}Z(0,t_{r+1})\Phi(t_{r+1},s)h(s,\varphi_1(s),\varphi_1(\gamma(s)))ds\}\notag
\end{align}
By Lemma 4.1, we have that
\begin{align}
\varphi_1(t)&= Z(t,0)(\varphi_1(0)+c_0) \notag\\
&+ \sum_{r=-\infty}^{j}\int_{t_r}^{\zeta_r}Z(t,0)PZ(0,t_r)\Phi(t_r,s)h(s,\varphi_1(s),\varphi_1(\gamma(s)))ds \notag\\
&+ \sum_{r=-\infty}^{j-1}\int_{\zeta_r}^{t_{r+1}}Z(t,0)PZ(0,t_{r+1})\Phi(t_{r+1},s)h(s,\varphi_1(s),\varphi_1(\gamma(s)))ds \notag\\
&- \sum_{r=j+1}^{+\infty}\int_{t_r}^{\zeta_r}Z(t,0)(I-P)Z(0,t_r)\Phi(t_r,s)h(s,\varphi_1(s),\varphi_1(\gamma(s)))ds \notag\\
&- \sum_{r=j}^{+\infty}\int_{\zeta_r}^{t_{r+1}}Z(t,0)(I-P)Z(0,t_{r+1})\Phi(t_{r+1},s)h(s,\varphi_1(s),\varphi_1(\gamma(s)))ds \notag\\
&+ \int_{\zeta_j}^{t}\Phi(t,s)h(s,\varphi_1(s),\varphi_1(\gamma(s)))ds \notag \\
&\triangleq Z(t,0)(\varphi_1(0)+c_0)+J, \notag
\end{align}
where
\begin{align}
c_0 &=\int_{0}^{\zeta_i}\Phi(0,s)h(s,\varphi_1(s),\varphi_1(\gamma(s)))ds  \notag\\
&- \sum_{r=-\infty}^{i}\int_{t_r}^{\zeta_r}PZ(0,t_r)\Phi(t_r,s)h(s,\varphi_1(s),\varphi_1(\gamma(s)))ds \notag\\
&- \sum_{r=-\infty}^{i-1}\int_{\zeta_r}^{t_{r+1}}PZ(0,t_{r+1})\Phi(t_{r+1},s)h(s,\varphi_1(s),\varphi_1(\gamma(s)))ds, \notag \\
&+\sum_{r=i+1}^{+\infty}\int_{t_r}^{\zeta_r}(I-P)Z(0,t_r)\Phi(t_r,s)h(s,\varphi_1(s),\varphi_1(\gamma(s)))ds \notag\\
&+ \sum_{r=i}^{+\infty}\int_{\zeta_r}^{t_{r+1}}(I-P)Z(0,t_{r+1})\Phi(t_{r+1},s)h(s,\varphi_1(s),\varphi_1(\gamma(s)))ds. \notag
\end{align}
Similar to the computation of  $|T\varphi(t)|$, we could prove that $J$ is bounded. Thus $Z(t,0)(\varphi_1(0)+c_0)$ is a bounded solution of system \eqref{syslz}. From Lemma 2.3, we have
$$\varphi_1(0)+c_0=0.$$
Thus
$$
\varphi_1(t)=\int_{-\infty}^{t}\tilde{G}_1(t,s)h(s,\varphi_1(s),\varphi_1(\gamma(s)))ds-\int_{t}^{+\infty}\tilde{G}_2(t,s)h(s,\varphi_1(s),\varphi_1(\gamma(s)))ds.
$$
Furthermore,
\begin{align}
&|\varphi_1(t)-\varphi_0(t)| \notag \\
&\leq|\int_{-\infty}^{t}\tilde{G}_1(t,s)[h(s,\varphi_1(s),\varphi_1(\gamma(s)))-h(s,\varphi_0(s),\varphi_0(\gamma(s)))]ds| \notag \\
&+|\int_{t}^{+\infty}\tilde{G}_2(t,s)[h(s,\varphi_1(s),\varphi_1(\gamma(s)))-h(s,\varphi_0(s),\varphi_0(\gamma(s)))]ds|  \notag\\
&\leq\int_{-\infty}^{t}K l\tilde{\rho}(M)e^{-\alpha(t-s)}(|\varphi_1(s)-\varphi_0(s)|+|\varphi_1(\gamma(s))-\varphi_0(\gamma(s))|)ds\notag\\
&+\int_{t}^{+\infty}Kl\tilde{\rho}(M)e^{\alpha(t-s)}(|\varphi_1(s)-\varphi_0(s)|+|\varphi_1(\gamma(s))-\varphi_0(\gamma(s))|)ds\notag\\
&\leq 4Kl\tilde{\rho}(M)\alpha^{-1}\|\varphi_1-\varphi_0\|  \notag\\
&\leq \frac{1}{2}\|\varphi_1-\varphi_0\|.\notag
\end{align}
Therefore
$$
 \|\varphi_1-\varphi_0\| \leq \frac{1}{2}\|\varphi_1-\varphi_0\|,
$$
which implies that $\varphi_1(t)=\varphi_0(t)$. This completes the proof.   $\Box$

\section{The preliminaries for the proof of Theorem 2}
\no

In this section, we give some preliminaries for the proof of Theorem 2.

\subsection{The solutions of subsystems }
\no

From Lemma 2.1, we have the following lemma.
\th{Lemma 5.1}
\emph{Assume that conditions \textbf{(A, $\mathfrak{B}$, $\mathfrak{C}$)} are fulfilled, then $J_k(t,s)(k=1,2)$ is nonsingular for any $t,s\in \bar{I}_r$ and the matrices
$Z_k(t,s)$ and $Z_k(t,s)^{-1}(k=1,2)$ are well defined for any $t, s\in\mathbb{R}$. If $t,s\in \bar{I}_r$, then
\begin{equation*}
|\Phi_1(t,s)|\leq\rho(A),\quad |\Phi_2(t,s)|\leq\rho(B),
\end{equation*}
\begin{equation*}
|Z_1(t,s)|\leq\rho_0(A), \quad  |Z_2(t,s)|\leq\rho_0(B),
\end{equation*}
where $\rho(\cdot)$ is defined in \eqref{defrho_AB} and $\rho_0(\cdot)$ is defined in \eqref{defalpha_0AB}.}

Lemma 5.1 ensures the continuity of solutions of subsystems \eqref{sysnlx} and  \eqref{sysnly} on $\mathbb{R}$. Moreover, we give the following remark.

\th{Remark 5.1} The fundamental matrix  $\Phi(t)$ of system
$\left(\begin{array}{c}x'(t)\\
  y'(t)\end{array}\right)=\left(\begin{array}{c}A(t)x(t)\\
  B(t)y(t)\end{array}\right) $
  with $\Phi(0)=I$ , and the transition matrix $Z(t,s)$ of system \eqref{syslxy}
  have the following form
$$\Phi(t,s)=\left(\begin{array}{c}\Phi_1(t,s)\\
 0\end{array}\begin{array}{c}0\\
 \Phi_2(t,s)\end{array}\right),\quad
Z(t,s)=\left(\begin{array}{c}Z_1(t,s)\\
 0\end{array}\begin{array}{c}0\\
 Z_2(t,s)\end{array}\right).$$

From Proposition 2.1, for any $t\in I_j$, $\tau\in I_i,$ the solution of subsystem \eqref{syslx} with $x(\tau)=\xi$ is defined on $\mathbb{R}$ and is given by
\begin{equation}\label{sollx}
x(t)=Z_1(t,\tau)\xi,
\end{equation}
and the solution of subsystem \eqref{sysly} with $y(\tau)=\eta$ can be represented as
\begin{equation}\label{solly}
y(t)=Z_2(t,\tau)\eta.
\end{equation}

From Proposition 2.2, for any $t\in I_j$, $\tau\in I_i$ and $t>\tau$, the solution of subsystem \eqref{sysnlx} with $x(\tau)=\xi$ is defined on $\mathbb{R}$ and is given by
\begin{eqnarray}\label{solnlx}
x(t)&=&Z_1(t,\tau)\xi+\int_\tau^{\zeta_i}Z_1(t,\tau)\Phi_1(\tau,s)(f(s)+\phi(s))ds
+\sum\limits_{r=i+1}^j\int_{t_r}^{\zeta_r}Z_1(t,t_r)\Phi_1(t_r,s)(f(s)+\phi(s))ds \nonumber\\
& & +\sum\limits_{r=i}^{j-1}\int_{\zeta_r}^{t_{r+1}}Z_1(t,t_{r+1})\Phi_1(t_{r+1},s)(f(s)+\phi(s))ds
+\int_{\zeta_j}^{t}\Phi_1(t,s)(f(s)+\phi(s))ds \nonumber\\
&\triangleq&Z_1(t,\tau)\xi+\int_\tau^{t}G_1(t,s)(f(s)+\phi(s))ds,
\end{eqnarray}
where $f(s)=f(s,x(s),x(\gamma(s)))$,  $\phi(s)=\phi(s,y(s),y(\gamma(s)))$ and
$$
G_1(t,s,\tau)= \left\{ \begin{array}{ccc}
Z_1(t,\tau)\Phi_1(\tau,s), & \textrm{ if } & s\in[\tau,\zeta_i] \textrm{ or } s\in[\zeta_i,\tau],\\
Z_1(t,t_r)\Phi_1(t_r,s), & \textrm{ if } & s\in[t_r,\zeta_r) \textrm{ for } r=i+1,\cdots,j,\\
Z_1(t,t_{r+1})\Phi_1(t_{r+1},s) & \textrm{ if } & s\in[\zeta_r,t_{r+1}) \textrm{ for } r=i,\cdots,j-1,\\
\Phi_1(t,s) & \textrm{ if } & s\in[\zeta_j,t] \textrm{ or } s\in[t,\zeta_j].
\end{array}
\right.
$$.

Similarly, if $t>\tau$, the solution of subsystem \eqref{sysnly} with $y(\tau)=\eta$ can be represented as
\begin{eqnarray}\label{solnly}
y(t)&=&Z_2(t,\tau)\eta+\int_\tau^{\zeta_i}Z_2(t,\tau)\Phi_2(\tau,s)(g(s)+\psi(s))ds
+\sum\limits_{r=i+1}^j\int_{t_r}^{\zeta_r}Z_2(t,t_r)\Phi_2(t_r,s)(g(s)+\psi(s))ds\nonumber\\
& & +\sum\limits_{r=i}^{j-1}\int_{\zeta_r}^{t_{r+1}}Z_2(t,t_{r+1})\Phi_2(t_{r+1},s)(g(s)+\psi(s)))ds
+\int_{\zeta_j}^{t}\Phi_2(t,s)(g(s)+\psi(s))ds\nonumber\\
&=&Z_2(t,\tau)\eta+\int_\tau^{t}G_2(t,s)(g(s)+\psi(s))ds,
\end{eqnarray}
where $g(s)=g(s,x(s),x(\gamma(s)))$,  $\psi(s)=\psi(s,y(s),y(\gamma(s)))$ and $G_2(t,s)$ can be defined in the same way as $G_1(t,s)$.

\th{Remark 5.2 \ } We could obtain $G_k(t,s)$ $(k=1,2)$ for $t<\tau$ by replacing $r=i+1,\cdots,j,$ and $r=i,\cdots,j-1,$ with $r=j+1,\cdots,i,$ and $r=j,\cdots,i-1$, in the definitions of $G_k(t,s)$ $(t>s, k=1,2)$, respectively.  From Remark 2.1, one could obtain the solution formulas of subsystems \eqref{sysnlx}  and \eqref{sysnly}   for the case $t<\tau$.

\subsection{Some lemmas}
\no
\th{Lemma \ 5.2}
\emph{If condition \textbf{($\mathfrak{D}$)} holds, for $t\in \mathbb{R}$ and  $s\in \mathbb{R}$, then
\begin{equation*}
|G_1(t,s)|\leq K\tilde{\rho}(A) e^{-\alpha(t-s)}, \quad t\geqslant  s, \qquad |G_2(t,s)|\leq K\tilde{\rho}(B) e^{\alpha(t-s)},\quad t<s,\\
\end{equation*}
where $\tilde{\rho}(\cdot)$ is defined in Theorem 1, $\alpha$ is  in \textbf{($\mathfrak{D}$)} and $\theta$ is  in \textbf{(A4)}.  }
\proof
We just prove the first inequality.

Suppose that $t\in I_j$, $\tau\in I_i$ and $t\geqslant s$.

{\bf Case $1$.}  $t\geqslant\tau$.

Without loss of generality, we assume that $t_i\leq\tau\leq\zeta_i\leq t_{i+1}\leq\cdots t_j\leq\zeta_j\leq t$.

If $s \in [\tau, \zeta_i]$, due to \textbf{(A4)}, we have $s-\tau\leq\theta$. It follows from \textbf{($\mathfrak{D}$)} and Lemma 5.1 that
$$|G_1(t,s)|=|Z_1(t,\tau)\Phi_1(\tau,s)|\leq Ke^{-\alpha(t-\tau)}\rho(A)\leq Ke^{-\alpha(t-s)}e^{\alpha\theta}\rho(A).$$

If $s \in [t_r,\zeta_r]$ ($r=i+1,\cdots,j$), then $s-t_r\leq\theta$. In view of \textbf{($\mathfrak{D}$)} and  Lemma 5.1, we have
$$
|G_1(t,s)|=|Z_1(t,t_r)\Phi_1(t_r,s)|\leq Ke^{-\alpha(t-t_r)}\rho(A)\leq Ke^{-\alpha(t-s)}e^{\alpha\theta}\rho(A).
$$

If $s \in [\zeta_r, t_{r+1}]$ ($r=i,\cdots,j-1$), similar to the above inequality, we have the same conclusion.

If $s \in [\zeta_j, t]$, owing to \textbf{(A4)}, we have $t-s\leq\theta$. It follows from Lemma 5.1 and $K\geqslant 1$ that
\begin{equation}\label{temp5.2-1}
|G_1(t,s)|=|\Phi_1(t,s)|\leq\rho(A)\leq Ke^{-\alpha(t-s)}e^{\alpha\theta}\rho(A).
\end{equation}

{\bf Case $2$.}  $t\leq\tau$.

By the definition of $G_1(t,s)$ we have $s\in [\min(t,\zeta_j), \max(\tau, \zeta_i)]$.

If $t\leq\zeta_j$, then $t<s$ which contradicts to our assumption that $t\geqslant s$. Thus, we just consider the case that $\zeta_j\leq t$. We divide the discussion into two subcases.

{\bf Subcase $2.1$.}  $\zeta_j\leq t\leq t_{j+1}\leq\tau$.

For $t\geqslant s$, the only possibility is that $s\in [\zeta_j,t]$. Similar to \eqref{temp5.2-1}, we have
$$
|G_1(t,s)|=|\Phi_1(t,s)|\leq Ke^{-\alpha(t-s)}e^{\alpha\theta}\rho(A).
$$

{\bf Subcase $2.2$.}  $\zeta_j\leq t\leq\tau\leq t_{j+1}$.

If $t\geqslant s$, then $s\in [\zeta_j, t]$ or $s \in [\zeta_j, \tau]$.

When $s\in [\zeta_j, t]$, similar to \eqref{temp5.2-1}, we get
$$
|G_1(t,s)|\leq Ke^{-\alpha(t-s)}e^{\alpha\theta}\rho(A).
$$

When $s\in [\zeta_j, \tau]$, we have $s\in \bar{I}_j$. Since $t\geqslant s$, following \textbf{($\mathfrak{D}$)} and Lemma 5.1, we obtain
$$|G_1(t,s)|=|Z_1(t,\tau)\Phi_1(\tau,s)|=|Z_1(t,s)Z_1(s,\tau)\Phi_1(\tau,s)|\leq Ke^{-\alpha(t-s)}\rho_0(A)\rho(A).$$

Note that $\tilde{\rho}(A)=\max(\rho(A)\rho_0(A), \rho(A)e^{\alpha\theta})$, we complete the proof.   $\Box$

Similar to Lemma 2.3, we have the following:
\th{Lemma \ 5.3 }\emph{Assume that condition \textbf{($\mathfrak{D}$)} holds, then
$$
\lim_{t\rightarrow -\infty}|Z_1(t,\tau)|=+\infty, \quad \lim_{t\rightarrow +\infty}|Z_2(t,\tau)|=+\infty, \quad \forall \tau\in \mathbb{R}.
$$
Moreover, the unique bounded solution in $\mathbb{R}$ of subsystem \eqref{syslx} (subsystem \eqref{sysly}) is trivial. }
\proof The proof is similar to that of Lemma 2.3  and so it is omitted.  $\Box$

\th{Lemma 5.4 (\cite{PintoRobledo15}, Lemma 5.1)} \emph{Let $t \mapsto z(t,\tau,\xi)$ and $t \mapsto z(t,\tau,\xi')$ be the solutions of system \eqref{sysnlxy-phi} passing respectively through $\xi$ and $\xi'$ at $t = \tau$. If \eqref{eq3.3} is valid, then it follows that
$$
|z(t,\tau,\xi')-z(t,\tau,\xi)| \leq |\xi-\xi'|e^{p(\ell)|t-\tau|}
$$
where $z(t,\cdot)=(x(t,\cdot), y(t,\cdot))^T$ and $p(\ell)$ is defined by
$$
p(\ell)= \eta_1 + \frac{\eta_2e^{\eta_1\theta}}{1-\upsilon} \quad with \quad \eta_1 = M + \ell, \quad \eta_2 = M_0 + \ell,
$$
and $\upsilon \in [0,1)$ is defined by \eqref{eq3.3}.  }

\th{Remark 5.3 \ } If $h(t,z(t),z(\gamma(t)))=0$, take $\ell=0$, Lemma 5.4 reduces to Lemma 5.2 in \cite{PintoRobledo15}.  Moreover, since  $p(\ell)>p(0)$ and $F(\ell,\theta)\geqslant F(0,\theta)$ in \eqref{eq3.3}, Lemma 5.4 is also valid for system \eqref{syslxy}.

\th{Lemma \ 5.5 (DEPCAG Gronwall inequality \cite{Pinto-MCM09,PintoJDEQ11, Coronel15})} \emph{ Let $\varrho,\eta: \mathbb{R} \rightarrow [0, \infty)$ be two functions such that $u$ is continuous and $¦Ç$ is locally integrable satisfying
$$
\bar{\theta} = \sup_{i\in\mathbb{Z}}\left\{\theta_i:\theta_i:=2\int_{I_i}\eta(s)ds\right\} < 1
$$
Suppose that for $\tau \leq t$ or $t \leq \tau$, we have the inequality
$$
\varrho(t) \leq \varrho(\tau) + \left|\int_{\tau}^{t}\eta(s)[\varrho(s)+\varrho(\gamma(s))]ds\right|.
$$
Then
$$
\varrho(t) \leq \varrho(\tau)\exp\left\{\tilde{\theta}\int_{\tau}^{t}\eta(s)ds\right\},
$$
$$
\varrho(\gamma(t)) \leq (1-\bar{\theta})^{-1}\varrho(\tau)\exp\left\{\tilde{\theta}\int_{\tau}^{t}\eta(s)ds\right\},
$$
where $\tilde{\theta}=\frac{2-\bar{\theta}}{1-\bar{\theta}}$.   }

\section{System \eqref{sysnlxy-fg} is topologically conjugate to system \eqref{syslxy}}

Suppose that
$
{\left(
\begin{array}{c}
    X(t,t_0,x_0) \\
    Y(t,t_0,x_0,y_0)
\end{array}
\right)}
$
is the solution of system \eqref{sysnlxy-fg} satisfying that
$
{\left(
\begin{array}{c}
    X(t_0) \\
    Y(t_0)
\end{array}
\right)} =
{\left(
\begin{array}{c}
    x_0 \\
    y_0
\end{array}
\right)}
$
and
$
{\left(
\begin{array}{c}
    u(t,t_0,\xi) \\
    v(t,t_0,\eta)
\end{array}
\right)}
$
is the solution of system \eqref{syslxy} satisfying that
$
{\left(
\begin{array}{c}
    u(t_0) \\
    v(t_0)
\end{array}
\right)} =
{\left(
\begin{array}{c}
    \xi \\
    \eta
\end{array}
\right)},
$
where $t_0\in\mathbb{R}$, $x_0,\xi\in\mathbb{R}^{n_1}$, $y_0, \eta\in\mathbb{R}^{n_2}.$

\th{Lemma 6.1} \emph{ For any $t\geq t_0$, the following inequalities hold:
$$
|X(t,t_0,x_0)|\leq |x_0|e^{-\alpha_0(t-t_0)},
$$
$$
|X(\gamma(t),t_0,x_0)|\leq (1-\bar{\theta})e^{\alpha_0\theta}|x_0|e^{-\alpha_0(t-t_0)},
$$
where $\alpha_0$ is defined in \eqref{eq3.2+2}.   }

\proof
From \eqref{sysnlx} we get
$$
X(t,t_0,x_0) =Z_1(t,t_0)x_0+\int_{t_0}^{t}G_1(t,s)f(s,X(s,t_0,x_0), X(\gamma(s),t_0,x_0))ds.
$$
It follows from condition \textbf($\mathfrak{D}$) and Lemma 5.2 that
$$
|X(t,t_0,x_0)|\leq e^{-\alpha(t-t_0)}|x_0|+l\tilde{\rho}(A)\int_{t_0}^{t}e^{-\alpha(t-s)}(|X(s)|+|X(\gamma(s))|)ds.
$$
Thus
\begin{align}
&\quad  e^{\alpha t}|X(t,t_0,x_0)|  \notag \\
&\leq e^{\alpha t_0}|x_0|+ l\tilde{\rho}(A)\int_{t_0}^{t}(e^{\alpha s}|X(s)|+e^{\alpha\theta}e^{\alpha \gamma(s)}|X(\gamma(s))|)ds \notag \\
&\leq e^{\alpha t_0}|x_0|+ l\tilde{\rho}(A)e^{\alpha\theta}\int_{t_0}^{t}(e^{\alpha s}|X(s)|+e^{\alpha \gamma(s)}|X(\gamma(s))|)ds. \notag
\end{align}
Applying Lemma 5.5 to $\varrho(t)=e^{\alpha t}|X(t,t_0,x_0)|$ and $\eta(t)=1$, we obtain that
$$
|X(t,t_0,x_0)| \leq |x_0|e^{-\alpha(t-t_0)+2 l\tilde{\rho}(A)e^{\alpha\theta}(t-t_0)}
$$
and
$$
|X(\gamma(t),t_0,x_0)| \leq (1-\bar{\theta})|x_0|e^{-\alpha(\gamma(t)-t_0)+2 l\tilde{\rho}(A)e^{\alpha\theta}(\gamma(t)-t_0)}.
$$
Thus
$$
|X(t,t_0,x_0)| \leq e^{-\alpha_0(t-t_0)}|x_0|,
$$
and
$$
|X(\gamma(t),t_0,x_0)|  \leq (1-\bar{\theta})e^{-\alpha_0(\gamma(t)-t_0)}|x_0| \leq (1-\bar{\theta})e^{\alpha_0\theta}e^{-\alpha_0(t-t_0)}|x_0|.\Box
$$

\th{Lemma 6.2} \emph{ For any fixed $t_0\in\mathbb{R}$,  $x_0, \xi\in\mathbb{R}^{n_1}$, there exists a unique $T(t_0,x_0)$ and $S(t_0,\xi)\in\mathbb{R}$, such that
$$
|X(T(t_0,x_0),t_0,x_0)|=1, \quad T(t_0,x_0) \rightarrow -\infty \quad \text{when} \quad x_0 \rightarrow 0,
$$
$$
|u(S(t_0,\xi),t_0,\xi|=1, \quad  S(t_0,\xi) \rightarrow -\infty  \quad \text{when} \quad  \xi \rightarrow 0.
$$   }

\proof
From Lemma 6.1, we have that $|X(t,t_0,x_0)| \leq |x_0|e^{-\alpha_0 (t-t_0)}$ when $t\geqslant t_0$, where $\alpha_0$ is defined in Lemma 6.1. If $x_0 \neq 0$ and $t \rightarrow +\infty,$ then
$$
\quad |X(t,t_0,x_0)| \rightarrow 0.
$$
If $t\geqslant\tau$,
\begin{equation}\label{Lemma4.1-temp1}
|X(t,t_0,x_0)|=|X(t,\tau,X(\tau,t_0,x_0))| \leq |X(\tau,t_0,x_0)|e^{-\alpha_0(t-\tau)}.
\end{equation}

Thus, for the fixed $t_0$ and $x_0$, $|X(t,t_0,x_0)|$ is a strictly monotonous decreasing function about $t$.
If $t$ is fixed and $\tau \rightarrow -\infty$, then
$$
e^{-\alpha_0(t-\tau)} \rightarrow 0.
$$
Thus
$$|X(\tau,t_0,x_0)| \rightarrow +\infty \quad \text{when} \quad \tau \rightarrow -\infty.$$

Therefore, there exists a unique time $T(t_0,x_0)$ such that $|X(T(t_0,x_0),t_0,x_0)|=1$.  Moreover, when $x_0 \rightarrow 0$, $T(t_0,x_0) \rightarrow -\infty.$

By condition \textbf{$\mathfrak{D}$}, for $t>t_0$, we have
$$|u(t,t_0,\xi)|=|Z(t,t_0)\xi| \leq e^{-\alpha (t-t_0)}|\xi|.$$
Thus when $t \rightarrow +\infty$,
$$|Z(t,t_0)\xi| \rightarrow 0.$$

Similar to  \eqref{Lemma4.1-temp1}, we could obtain that for fixed $t_0$ and $\xi$, $|Z(t,t_0)\xi|$ is a strictly  monotonous decreasing function about $t$. Moreover, when $t \rightarrow -\infty$,
$$|Z(t,t_0)\xi| \rightarrow +\infty.$$
Therefore, for a fixed $\xi\in\mathbb{R}^{n_1},\xi \neq 0$, there exists a unique time $S(t_0,\xi)$ such that
$$|Z(S(t_0,\xi),t_0,\xi|=1,$$
and
$$
S(t_0,\xi) \rightarrow -\infty \quad  \text{when} \quad \xi \rightarrow 0.   \quad  \Box
$$

\th{Lemma 6.3} \emph{ For any $x_0 \neq 0$,  $\xi \neq 0$ and $t\in\mathbb{R}$, we have
$$
T(t,X(t,t_0,x_0)) = T(t_0,x_0),
$$
$$
S(t,u(t,\tau,\xi)) = S(\tau,\xi).
$$   }

\proof
It follows from Lemma 6.2 that
$$
1=|X(T(t,X(t,t_0,x_0)),t,X(t,t_0,x_0))|=|X(T(t,X(t,t_0,x_0)),t_0,x_0)|.
$$
From $|X(T(t_0,x_0),t_0,x_0))|=1$ and Lemma 6.2, we get
$$
T(t,X(t,t_0,x_0)) = T(t_0,x_0).
$$
The second equality can be proved in a similar way.   $\Box$

\th{Lemma 6.4}   \emph{
For any $t_0\in\mathbb{R}$, $x_0\in\mathbb{R}^{n_1}$, the following inequality holds.
$$|\int_{t_0}^{+\infty}G_2(t_0,s)g(s,X(s,t_0,x_0),X(\gamma(s),t_0,x_0))ds|\leq K\lambda \tilde{\rho}(B)((\alpha+\rho_0)^{-1}+\alpha^{-1})|x_0|.$$    }

\proof
From Lemmas 5.2 and 6.1, we get
\begin{align}
 &|\int_{t_0}^{+\infty}G_2(t_0,s)g(s,X(s,t_0,x_0),X(\gamma(s),t_0,x_0))ds| \notag \\
& \leq \int_{t_0}^{+\infty}K\lambda\tilde{\rho}(B)e^{\alpha (t_0-s)}(|X(s,t_0,x_0)|+|X(\gamma(s),t_0,x_0)|) ds \notag \\
&\leq \int_{t_0}^{+\infty}K\lambda\tilde{\rho}(B)e^{\alpha (t_0-s)}e^{-\alpha_0(s-t_0)}(|x_0|+ (1-\bar{\theta})e^{\alpha_0\theta}|x_0|) ds \notag \\
&\leq K\lambda(\alpha+\alpha_0)^{-1}\tilde{\rho}(B)(1+(1-\bar{\theta})e^{\alpha_0\theta})|x_0|.  \notag
\end{align}

\th{Definition 6.2}
For any $t \in \mathbb{R}$, $\xi \in \mathbb{R}^{n_1}$ and $\eta \in \mathbb{R}^{n_2}$, we define
$L_1:\mathbb{R}\times\mathbb{R}^{n_1} \rightarrow \mathbb{R}^{n_1}$,  $L_2:\mathbb{R}\times\mathbb{R}^{n_1}\times\mathbb{R}^{n_2} \rightarrow \mathbb{R}^{n_2}$  and  $L:\mathbb{R}\times\mathbb{R}^{n_1}\times\mathbb{R}^{n_2} \rightarrow \mathbb{R}^{n}$ as follows:
$$
L_1(t,\xi)= \left\{
\begin{array}{c}
    X(t,S(t,\xi),u(S(t,\xi),t,\xi)) \quad  \xi\neq 0,  \\
    0 \qquad \qquad   \qquad   \qquad   \qquad   \qquad \xi=0,
\end{array}
\right.
$$
$$
L_2(t,\xi,\eta) = \eta - \int_{t}^{+\infty}G_2(t,s)g\Big{(}s,X(s,t,L_1(t,\xi)),X(\gamma(s),t,L_1(t,\xi))\Big{)}ds,
$$
and
$$
L(t,\xi,\eta)=\
{\left(
\begin{array}{c}
    L_1(t,\xi) \\
    L_2(t,\xi,\eta)
\end{array}
\right)}.\
$$

\th{Lemma 6.5}
\emph{ $L_1(t,\xi)$ is a continuous function of $\xi$ and $L_1(t,u(t,\tau,\xi))=X(t,\tau,L_1(\tau,\xi)).$   }

\proof
By Lemma 6.2, we have
$$S(t,\xi) \rightarrow -\infty    \quad    \text{when}   \quad    \xi \rightarrow 0.$$
When $\xi \rightarrow 0$, it follows from Lemma 6.1 that
$$|X(t,S(t,\xi),u(S(t,\xi),t,\xi))| \leq |u(S(t,\xi),t,\xi))|e^{-\alpha_0(t-S(t,\xi))}=e^{-\alpha_0(t-S(t,\xi))} \rightarrow 0.$$
Hence, $L_1(t,\xi)$ is  a continuous function of $\xi$.

Furthermore, from Lemma 6.3, we have that
\begin{align}
L_1(t,u(t,\tau,\xi)) \
&=X(t,S(t,u(t,\tau,\xi)),u(S(t,u(t,\tau,\xi)),t,u(t,\tau,\xi))) \notag\\
&=X(t,S(\tau,\xi),u(S(\tau,\xi),\tau,\xi))    \notag\\
&=X(t,\tau,X(\tau,S(\tau,\xi),u(S(\tau,\xi),\tau,\xi))) \notag\\
&=X(t,\tau,L_1(\tau,\xi)).       \quad \Box         \notag
\end{align}

\th{Lemma 6.6}  \emph{
$
{\left(
\begin{array}{c}
    L_1(t,u(t,\tau,\xi)) \\
    L_2(t,u(t,\tau,\xi),v(t,\tau,\eta))
\end{array}
\right)} = \
{\left(
\begin{array}{c}
    X(t,\tau,L_1(\tau,\xi)) \\
    Y(t,\tau,L_1(\tau,\xi),L_2(\tau,\xi,\eta))
\end{array}
\right).}
$          }

\proof
Due to Lemma 6.5, we get
$$L_1(t,u(t,\tau,\xi))=X(t,\tau,L_1(\tau,\xi)).$$

\begin{align}\label{le6.5-temp}
&L_2(t,u(t,\tau,\xi),v(t,\tau,\eta)) \notag\\
&= v(t,\tau,\eta) - \int_{t}^{+\infty}G_2(t,s)g\Big{(}s,X(s,t,L_1(t,u(t,\tau,\xi))),X(\gamma(s),t,L_1(t,u(t,\tau,\xi)))\Big{)}ds \notag\\
&= v(t,\tau,\eta) - \int_{t}^{+\infty}G_2(t,s)g\Big{(}s,X(s,t,X(t,\tau,L_1(\tau,\xi))),X(\gamma(s),t,X(t,\tau,L_1(\tau,\xi)))\Big{)}ds \notag\\
&= v(t,\tau,\eta) - \int_{t}^{+\infty}G_2(t,s)g\Big{(}s,X(s,\tau,L_1(\tau,\xi)),X(\gamma(s),\tau,L_1(\tau,\xi))\Big{)}ds.
\end{align}

Denote $J(t)=- \int_{t}^{+\infty}G_2(t,s)g\Big{(}s,X(s,\tau,L_1(\tau,\xi)),X(\gamma(s),\tau,L_1(\tau,\xi))\Big{)}ds.$
Suppose $t\in I_j$, we obtain
\begin{align}
J'(t)&= - B(t)\int_{t}^{+\infty}G_2(t,s)g\Big{(}s,X(s,\tau,L_1(\tau,\xi)),X(\gamma(s),\tau,L_1(\tau,\xi))\Big{)}ds \notag\\
& - B_0(t)\int_{\gamma(t)}^{+\infty}G_2(t,s)g\Big{(}s,X(s,\tau,L_1(\tau,\xi)),X(\gamma(s),\tau,L_1(\tau,\xi))\Big{)}ds  \notag\\
&+ g(t,X(t,\tau,L_1(\tau,\xi)),X(\gamma(t),\tau,L_1(\tau,\xi))). \notag
\end{align}

Furthermore, from \eqref{le6.5-temp}, we have
\begin{align}
&L_2'(t,u(t,\tau,\xi),v(t,\tau,\eta))\notag\\
=& B(t)L_2(t,u(t,\tau,\xi),v(t,\tau,\eta)) + B_0(t)L_2(\gamma(t),u(t,\tau,\xi),v(t,\tau,\eta)) \notag  \\
&+ g(t,X(t,\tau,L_1(\tau,\xi)),X(\gamma(t),\tau,L_1(\tau,\xi))).  \notag
\end{align}

Thus
$
{\left(
\begin{array}{c}
    L_1(t,u(t,\tau,\xi),v(t,\tau,\eta)) \\
    L_2(t,u(t,\tau,\xi),v(t,\tau,\eta))
\end{array}
\right)}
$
is a solution of system \eqref{sysnlxy-fg}.

From
\begin{align}
&{\left(
\begin{array}{c}
    L_1(t,u(t,\tau,\xi),v(t,\tau,\eta)) \\
    L_2(t,u(t,\tau,\xi),v(t,\tau,\eta))
\end{array}
\right)}|_{t=\tau}
= \
{\left(
\begin{array}{c}
    L_1(\tau,\xi) \\
    L_2(\tau,\xi,\eta)
\end{array}
\right)}   \notag
\end{align}
and
$$
{\left(
\begin{array}{c}
    X(t,\tau,L_1(\tau,\xi)) \\
    Y(t,\tau,L_1(\tau,\xi),L_2(\tau,\xi,\eta))
\end{array}
\right)}|_{t=\tau} = \
{\left(
\begin{array}{c}
    L_1(\tau,\xi) \\
    L_2(\tau,\xi,\eta))
\end{array}
\right),}
$$
we get the conclusion of the lemma.   $\Box$

\th{Definition 6.3}
For any $t\in \mathbb{R}$,  $x\in \mathbb{R}^{n_1}$  and  $y \in \mathbb{R}^{n_2}$, we denote
$
H(t,x,y) = \
{\left(
\begin{array}{c}
    H_1(t,x) \\
    H_2(t,x,y)
\end{array}
\right),}
$
where $H_1(t,x)$ and $H_2(t,x,y)$ are defined as
$$
H_1(t,x) = \
{\left\{
\begin{array}{c}
    u(t,T(t,x),X(T(t,x),t,x))   \quad x \neq 0, \\
    0 \qquad \qquad \qquad  \qquad   \qquad  \qquad x = 0,
\end{array}
\right.}
$$
and
$$
H_2(t,x,y) = y + \int_{t}^{+\infty}G_2(t,s)g(s,X(s,t,x),X(\gamma(s),t,x))ds.
$$

\th{Lemma 6.7} \emph{ $H_1(t,x)$ is a continuous function of $x$.    }
\proof
From \eqref{sollx}, we get
$$u(t,T(t,x),X(T(t,x),t,x))=Z_1(t,T(t,x))X(T(t,x),t,x),$$
which together with condition \textbf{$(\mathfrak{D})$} implies that
$$|u(t,T(t,x),X(T(t,x),t,x))|\leq e^{-\alpha(t-T(t,x))}|X(T(t,x),t,x)|\leq e^{-\alpha(t-T(t,x))}, \quad t\geqslant T(t,x).$$

From Lemma 6.2, we have that
$$T(t,x)\rightarrow -\infty, \quad \text{when} \quad x \rightarrow 0.$$
Thus $H_1(t,x)$ is a continuous function of $x$.    $\Box$

\th{Lemma 6.8}  \emph{
$$
{\left(
\begin{array}{c}
    H_1(t,X(t,t_0,x_0)) \\
    H_2(t,X(t,t_0,x_0),Y(t,t_0,x_0,y_0))
\end{array}
\right)} = \
{\left(
\begin{array}{c}
    u(t,t_0,H_1(t_0,x_0)) \\
    v(t,t_0,H_2(t,x_0,y_0))
\end{array}
\right).}
$$   }

\proof
From Lemma 6.3, we have
\begin{align}
H_1(t,X(t,t_0,x_0)) \
&= u(t, T(t, X(t,t_0,x_0)), X(T(t,X(t,t_0,x_0)), t, X(t,t_0,x_0))) \notag\\
&= u(t, T(t_0,x_0), X(T(t_0,x_0),t_0,x_0)) \notag\\
&= u(t, t_0, u(t_0, T(t_0,x_0), X(T(t_0,x_0), t_0, x_0))) \notag\\
&= u(t, t_0, H_1(t_0,x_0)). \notag
\end{align}

\begin{align}
&H_2(t, X(t,t_0,x_0), Y(t,t_0,x_0,y_0)) \notag  \\
&= Y(t,t_0,x_0,y_0) + \int_{t}^{+\infty}G_2(t,s)g\Big{(}s, X(s, t, X(t,t_0,x_0)), X(\gamma(s), t, X(t,t_0,x_0))\Big{)}ds \notag\\
&= Y(t,t_0,x_0,y_0) + \int_{t}^{+\infty}G_2(t,s)g(s, X(s,t_0,x_0), X(\gamma(s),t_0,x_0))ds \notag
\end{align}

Since
\begin{align}
&H_2'(t, X(t,t_0,x_0)) \notag \\
=& B(t)H_2(t, X(t,t_0,x_0), Y(t,t_0,x_0,y_0)) + B_0(t)H_2(\gamma(t), X(\gamma(t),t_0,x_0), Y(\gamma(t),t_0,x_0,y_0)), \notag
\end{align}
$H_2(t, X(t,t_0,x_0), Y(t,t_0,x_0,y_0))$ is a solution of system \eqref{sysly}.

Moreover,
$$
H_2(t, X(t,t_0,x_0), Y(t,t_0,x_0,y_0))|_{t=t_0} = H_2(t_0,x_0,y_0).
$$

Thus $H_2(t, X(t,t_0,x_0), Y(t,t_0,x_0,y_0))=v(t,t_0,H_2(t_0,x_0,y_0)).$ $\Box$

\th{Lemma 6.9} \emph{ For any $t_0 \in \mathbb{R}$,  $x_0 \in \mathbb{R}^{n_1}$,  $\tau \in \mathbb{R}$ and $\xi \in \mathbb{R}^{n_1}$, we have
$$
S(t_0, H_1(t_0,x_0)) = T(t_0,x_0), \quad  T(\tau, L_1(\tau,\xi)) = S(\tau,\xi).
$$        }

\proof
From the definition of $H_1$, we have
\begin{align}
&1=|u(S(t_0, H_1(t_0,x_0)), t_0, H_1(t_0,x_0))| \notag \\
&=|u(S(t_0, H_1(t_0,x_0)), t_0, u(t_0, T(t_0,x_0), X(T(t_0,x_0), t_0, x_0))) |  \notag \\
&=|u(S(t_0, H_1(t_0,x_0)), T(t_0,x_0), X(T(t_0,x_0), t_0, x_0))|,  \notag
\end{align}
which implies that
$$
S(t_0, H_1(t_0,x_0)) =S(T(t_0,x_0), X(T(t_0,x_0), t_0, x_0)).
$$
From
$$|u(T(t_0,x_0), T(t_0,x_0), X(T(t_0,x_0), t_0, x_0))| = |X(T(t_0,x_0), t_0, x_0)| = 1,$$
we obtain that
$$
S(T(t_0,x_0), X(T(t_0,x_0), t_0, x_0))= T(t_0,x_0).
$$
Thus
$$
S(t_0, H_1(t_0,x_0)) = T(t_0,x_0).
$$
Similarly, we could prove that $T(\tau, L_1(\tau,\xi)) = S(\tau,\xi).$   $\Box$

\th{Lemma 6.10} \emph{ For any $t_0 \in \mathbb{R}$,  $x_0 \in \mathbb{R}^{n_1}$ and $y_0 \in \mathbb{R}^{n_2}$,  we have
$$L(t_0,H(t_0,x_0))=(x_0,y_0)^T.$$  }

\proof
If $x_0=0$, it is easy to see that $L_1(t_0, H_1(t_0,x_0)) = x_0.$

If $x_0 \neq 0$, from Lemma 6.9 and the definitions of $L1$ and $H_1$, we get
\begin{align}
L_1(t_0, H_1(t_0,x_0)) \
&= X\Big{(}t_0, S(t_0, H_1(t_0,x_0)), u\big{(}S(t_0, H_1(t_0,x_0)), t_0, H_1(t_0,x_0)\big{)}\Big{)} \notag\\
&= X\Big{(}t_0, T(t_0,x_0), u\big{(}T(t_0,x_0), t_0, u(t_0, T(t_0,x_0), X(T(t_0,x_0), t_0, x_0))\big{)}\Big{)}  \notag\\
&= X\Big{(}t_0, T(t_0,x_0), u\big{(}T(t_0,x_0), T(t_0,x_0), X(T(t_0,x_0), t_0, x_0)\big{)}\Big{)}  \notag\\
&= X(t_0, T(t_0,x_0), X(T(t_0,x_0), t_0, x_0)) \notag\\
&= x_0, \notag
\end{align}
which together with the definitions of $L_2$ and $H_2$ implies that
\begin{align}
&\quad L_2(t_0, H_1(t_0,x_0), H_2(t_0,x_0,y_0)) \notag\\
&= H_2(t_0,x_0,y_0) - \int_{t_0}^{+\infty}G_2(t_0,s)g\Big{(}s, X(s, t_0, L_1(t_0, H_1(t_0,x_0))),X(\gamma(s), t_0, L_1(t_0, H_1(t_0,x_0)))\Big{)}ds \notag\\
&= y_0 + \int_{t_0}^{+\infty}G_2(t_0,s)g(s, X(s,t_0,x_0), X(\gamma(s),t_0,x_0))ds  \notag\\
&\qquad   - \int_{t_0}^{+\infty}G_2(t_0,s)g(s,X(s,t_0,x_0),X(\gamma(s),t_0,x_0))ds \notag\\
&= y_0. \notag   \quad \Box
\end{align}

\th{Lemma 6.11} \emph{ For any  $\tau \in \mathbb{R}$,  $\xi \in \mathbb{R}^{n_1}$  and  $\eta \in \mathbb{R}^{n_2}$, we have $$H(\tau,L(\tau,\xi,\eta))=(\xi,\eta)^T.$$     }

\proof
If $\xi=0$, it is obvious that $H_1(\tau,L_1(\tau,\xi))=\xi.$

If $\xi \neq 0$, by Lemma 6.9 and the definitions of $H_1$ and $L_1$, we obtain
\begin{align}
H_1(\tau,L_1(\tau,\xi)) \
&= u\Big{(}\tau, T(\tau,L_1(\tau,\xi)), X\big{(}T(\tau,L_1(\tau,\xi)), \tau, L_1(\tau,\xi)\big{)}\Big{)} \notag\\
&= u\Big{(}\tau, S(\tau,\xi), X\big{(}S(\tau,\xi), \tau, L_1(\tau,\xi)\big{)}\Big{)} \notag\\
&= u\Big{(}\tau, S(\tau,\xi), X\big{(}S(\tau,\xi), \tau, X(\tau, S(\tau,\xi), u(S(\tau,\xi), \tau, \xi))\big{)}\Big{)} \notag\\
&= u(\tau, S(\tau,\xi), u(S(\tau,\xi), \tau, \xi)) \notag\\
&= \xi. \notag
\end{align}

In what follows, we prove that $H_2(\tau,L_1(\tau,\xi),L_2(\tau,\xi,\eta))=\eta.$

For any $t \in \mathbb{R}, x \in \mathbb{R}^{n_1}$ and $y \in \mathbb{R}^{n_2}$,  due to Lemma 6.4, we have
\begin{align}
|H_2(t,x,y)-y| \
&= |\int_{t}^{+\infty}G_2(t,s)g(s,X(s,t,x),X(\gamma(s),t,x))ds| \notag\\
&\leq K\lambda \tilde{\rho}(B)((\alpha+\rho_0)^{-1}+\alpha^{-1})|x_0|. \notag
\end{align}

From Lemma 6.4 and the definition of $L_2$, we obtain
\begin{align}
|L_2(t,\xi,\eta)-\eta| \
&\leq |\int_{t}^{+\infty}G_2(t,s)g\Big{(}s,X(s,t,L_1(t,\xi)),X(\gamma(s),t,L_1(t,\xi))\Big{)}ds| \notag\\
&\leq K\lambda \tilde{\rho}(B)((\alpha+\rho_0)^{-1}+\alpha^{-1})|L_1(t,\xi)|. \notag
\end{align}

Thus, by Lemma 6.6  we get
\begin{align}
&J\triangleq|H_2(t, L_1(t,u(t,\tau,\xi)), L_2(t, u(t,\tau,\xi), v(t,\tau,\eta))) - v(t,\tau,\eta)| \notag\\
&\leq |H_2(t, L_1(t,u(t,\tau,\xi)), L_2(t, u(t,\tau,\xi), v(t,\tau,\eta))) - L_2(t, u(t,\tau,\xi), v(t,\tau,\eta))| \notag\\
&\quad + |L_2(t, u(t,\tau,\xi), v(t,\tau,\eta)) - v(t,\tau,\eta)| \notag\\
&\leq 2K\lambda \tilde{\rho}(B)((\alpha+\rho_0)^{-1}+\alpha^{-1})|L_1(t,u(t,\tau,\xi))|\notag\\
&\leq 2K\lambda \tilde{\rho}(B)((\alpha+\rho_0)^{-1}+\alpha^{-1})|X(t,\tau,L_1(\tau,\xi))|. \notag
\end{align}

It follows from Lemma 6.1 that
\begin{equation}\label{Lemma6.11-temp1}
J\leq 2K\lambda \tilde{\rho}(B)((\alpha+\rho_0)^{-1}+\alpha^{-1})|L_1(\tau,\xi)|e^{-\alpha_0(t-\tau)}, \quad t\geqslant\tau.
\end{equation}

From \eqref{solly},  Lemmas 6.6 and 6.8, we have
\begin{align}
&\quad H_2(t, L_1(t,u(t,\tau,\xi)), L_2(t,u(t,\tau,\xi),v(t,\tau,\eta)))  \notag\\
&= H_2(t, X(t,\tau,L_1(\tau,\xi)), Y(t,\tau,L_1(\tau,\xi)),L_2(\tau,\xi,\eta))) \notag\\
&= v(t, \tau, H_2(t,L_1(\tau,\xi),L_2(\tau,\xi,\eta))) \notag \\
&= Z_2(t,\tau)H_2(\tau, L_1(\tau,\xi), L_2(\tau,\xi,\eta)).   \notag
\end{align}

By \eqref{Lemma6.11-temp1} and $v(t,\tau,\eta)=Z_2(t,\tau)\eta$, we get
\begin{align}
&\quad|Z_2(t,\tau)\cdot\Big{(}H_2(\tau,L_1(\tau,\xi),L_2(\tau,\xi,\eta))-\eta\Big{)}| \notag\\
&= |H_2(t, L_1(t,u(t,\tau,\xi)), L_2(t, u(t,\tau,\xi), v(t,\tau,\eta))) - v(t,\tau,\eta)| \notag\\
&\leq 2K\lambda \tilde{\rho}(B)((\alpha+\rho_0)^{-1}+\alpha^{-1})|L_1(\tau,\xi)|e^{-\alpha_0(t-\tau)}, \quad t\geqslant\tau.   \notag
\end{align}

For  fixed $\tau$ and $\xi$, $L_1(\tau,\xi)$ is a fixed value. Thus the above equality is bounded when $t\geq \tau$. Moreover, it follows from  condition \textbf{$(\mathfrak{D})$} that the above equality is bounded when $t\leq \tau$.  Therefore, $Z_2(t,\tau)\cdot\Big{(}H_2(\tau,L_1(\tau,\xi),L_2(\tau,\xi,\eta))-\eta\Big{)}$ is a bounded solution of system \eqref{sysly}.

Since system \eqref{sysly} has an $\alpha$-exponential dichotomy, for fixed $\tau$, $\xi$ and $\eta$, it has a unique bounded solution, zero solution. Thus
$$H_2(\tau, L_1(\tau,\xi), L_2(\tau,\xi,\eta)) - \eta = 0.$$
That is   $H_2(\tau, L_2(\tau,\xi,\eta)) = \eta.$    $\Box$

\th{Lemma 6.12}  \emph{ System \eqref{sysnlxy-fg} is topologically conjugate to system \eqref{syslxy}.  }
\proof
It follows from Lemmas 6.10 and 6.11 that for a fixed $t$, $H(t,x,y): \mathbb{R}^{n_1}\times \mathbb{R}^{n_2} \rightarrow \mathbb{R}^{n}$ is a bijection and $H^{-1}(t,x,y)=L(t,x,y).$

According to Lemma 5.4 and Remark 5.3,  solutions of systems \eqref{sysnlxy-fg} and  \eqref{syslxy} are continuous with respect to initial values.

By the  definitions of  $H(t,\cdot)$ and $L(t,\cdot)$, and lemmas 6.5 and 6.7, we get that both $H(t,\cdot)$ and $L(t,\cdot)$ are continuous. Thus $H(t,\cdot)$ and $L(t,\cdot)$ are homeomorphisms of $\mathbb{R}^n$.

Moreover, Lemmas 6.6 and 6.8 imply that $H(t,\cdot)$ sends the solutions of system \eqref{sysnlxy-fg} onto those
of system \eqref{syslxy} and $L(t,\cdot)$ sends the solutions of system \eqref{syslxy} onto those of system \eqref{sysnlxy-fg}. Therefore, system \eqref{sysnlxy-fg} and system \eqref{syslxy} are topologically
conjugated. $\Box$

\section{System \eqref{sysnlxy-phi}  is topologically conjugate to system \eqref{sysnlxy-fg}}

First we introduce a new system
\begin{equation}\label{sysnlxy-pq}
\left\{\begin{array}{lc}
x' = A(t)x(t) + A_0(t)x(\gamma(t)) +f(t,x(t),x(\gamma(t))) + p(t,y(t),y(\gamma(t))) \\
y' = B(t)y(t) + B_0(t)x(\gamma(t))+ g(t,x(t),x(\gamma(t))) + q(t,y(t),y(\gamma(t))),
\end{array}\right.
\end{equation}
where $f(t,\cdot)$ and $g(t,\cdot)$ are defined in system \eqref{sysnlxy-phi}, $p:\mathbb{R}\times\mathbb{R}^{n_2}\times\mathbb{R}^{n_2}\rightarrow\mathbb{R}^{n_1}$ and $q:\mathbb{R}\times\mathbb{R}^{n_2}\times\mathbb{R}^{n_2}\rightarrow\mathbb{R}^{n_2}$ satisfying that
 for the $\delta$ and $\omega$ in \textbf{($\mathfrak{B}_2$)} and any $ t \in \mathbb{R}$, $y_1, y_2, \bar{y}_1,\bar{y}_2 \in \mathbb{R}^{n_2}$ such that
$$
|p(t,y_1, y_2)| \leq \delta, \quad |q(t,,y_1, y_2)| \leq \delta,
$$
$$
|p(t,y_1, y_2)-p(t,\bar{y}_1,\bar{y}_2)| \leq \omega(|y_1-\bar{y}_1|+|y_2-\bar{y}_2|),$$

$$|q(t,y_1,y_2)-q(t,\bar{y}_1,\bar{y}_2))| \leq \omega(|y_1-\bar{y}_1|+|y_2-\bar{y}_2|).
$$

\th{Lemma 7.1}  \emph{If \eqref{eq3.2+1} holds, then  there exists a unique function $\bar{H}(t,x,y): \mathbb{R} \times \mathbb{R}^{n_1+n_2} \rightarrow \mathbb{R}^{n_1+n_2}$ satisfying that }
\begin{description}
\item[(i)]
\emph{There exists a  constant $\bar{\sigma}>0$ such that
$$
|\bar{H}(t,x,y)-
(x, y)^T|
\leq   \bar{\sigma}.
$$
}
\item[(ii)]
\emph{ If
$
{\left(
\begin{array}{c}
    x(t) \\
    y(t)
\end{array}
\right)}
$
is a solution of system \eqref{sysnlxy-phi}, then $\bar{H}(t,x(t),y(t))$ is a solution of system \eqref{sysnlxy-pq}.
 }
\end{description}

\proof
For any fixed $\tau \in \mathbb{R}$, $\xi \in \mathbb{R}^{n_1}$  and $\eta \in \mathbb{R}^{n_2}$, suppose that
$
{\left(
\begin{array}{c}
    x(t,\tau,\xi,\eta) \\
    y(t,\tau,\xi,\eta)
\end{array}
\right)}
$
is a solution of system \eqref{sysnlxy-phi} satisfying
$
{\left(
\begin{array}{c}
    x(\tau,\tau,\xi,\eta) \\
    y(\tau,\tau,\xi,\eta)
\end{array}
\right)} =
{\left(
\begin{array}{c}
    \xi \\
    \eta
\end{array}
\right).}
$

Denote
$
z(t)=
{\left(
\begin{array}{c}
    z_1(t) \\
    z_2(t)
\end{array}
\right)}
$ where $z_1(t) \in \mathbb{R}^{n_1}$ and $z_2(t) \in \mathbb{R}^{n_2},$ $
W(t)=
{\left[
\begin{array}{cc}
    A(t) &  \\
     & B(t)
\end{array}
\right],}
$
$
W_0(t)=
{\left[
\begin{array}{cc}
    A_0(t) &  \\
     & B_0(t)
\end{array}
\right]}
$
 and
\begin{align}
&\quad \quad\bar{ h}(t,z(t),z(\gamma(t)),(\tau,\xi,\eta))\notag\\
&= {\left(
\begin{array}{c}
    \bar{h}_1(t,z(t),z(\gamma(t)),(\tau,\xi,\eta)) \\
    \bar{h}_2(t,z(t),z(\gamma(t)),(\tau,\xi,\eta))
\end{array}
\right)} \notag\\
&= {\left(
\begin{array}{c}
    f(t,x(t,\tau,\xi,\eta)+z_1(t),x(\gamma(t),\tau,\xi,\eta)+z_1(\gamma(t)))  \\
    g(t,x(t,\tau,\xi,\eta)+z_1(t),x(\gamma(t),\tau,\xi,\eta)+z_1(\gamma(t)))
\end{array}
\right)} \notag \\
&\quad\quad+ {\left(
\begin{array}{c}
    p(t,y(t,\tau,\xi,\eta)+z_2(t),y(\gamma(t),\tau,\xi,\eta)+z_2(\gamma(t)))  \\
    q(t,y(t,\tau,\xi,\eta)+z_2(t),y(\gamma(t),\tau,\xi,\eta)+z_2(\gamma(t)))
\end{array}
\right)} \notag \\
&\quad\quad+ {\left(
\begin{array}{c}
     - f(t,x(t,\tau,\xi,\eta),x(\gamma(t),\tau,\xi,\eta)) - \phi(t,y(t,\tau,\xi,\eta),y(\gamma(t),\tau,\xi,\eta)) \\
     - g(t,x(t,\tau,\xi,\eta),x(\gamma(t),\tau,\xi,\eta)) - \psi(t,y(t,\tau,\xi,\eta),y(\gamma(t),\tau,\xi,\eta))
\end{array}
\right).} \notag
\end{align}

From
$$|\bar{h}(t,z(t),z(\gamma(t)),(\tau,\xi,\eta))| \leq 2\lambda|z(t)| + 2\lambda|z(\gamma(t))|+4\delta, $$

\begin{align}
&\quad |\bar{h}(t,z(t),z(\gamma(t)),(\tau,\xi,\eta)) - \bar{h}(t,\bar{z}(t),\bar{z}(\gamma(t)),(\tau,\xi,\eta))| \notag\\
&\leq 2\omega |z(t)-\bar{z}(t)|+2\omega |z(\gamma(t))-\bar{z}(\gamma(t))|, \notag
\end{align}
and Theorem 1, we get that system
\begin{equation}\label{Lemma7.1-temp1}
z'(t)=
W(t)z(t) +W_0(t)z(\gamma(t))+ \bar{h}(t,z(t),z(\gamma(t)),(\tau,\xi,\eta))
\end{equation}
has a unique bounded solution for fixed $\tau$, $\xi$  and $\eta$.
We denote  by
$$
\chi(t,(\tau,\xi,\eta))
 = {\left(
\begin{array}{c}
    \chi_1(t,(\tau,\xi,\eta)) \\
    \chi_2(t,(\tau,\xi,\eta))
\end{array}
\right)}\quad \text{and} \quad |\chi(t,(\tau,\xi,\eta))| \leq \bar{\sigma},
$$
where $\chi_1(t,(\tau,\xi,\eta))\in \mathbb{R}^{n_1}$ and $\chi_2(t,(\tau,\xi,\eta))\in \mathbb{R}^{n_2}.$

For any  $t \in \mathbb{R}$,  $\xi \in \mathbb{R}^{n_1}$ and $\eta \in \mathbb{R}^{n_2}$, define
$$
\bar{H}(t,\xi,\eta) ={\left(
\begin{array}{c}
    \bar{H}_1(t,\xi,\eta) \\
    \bar{H}_2(t,\xi,\eta)
\end{array}
\right)}=
{\left(
\begin{array}{c}
    \xi + \chi_1(t,(t,\xi,\eta)) \\
    \eta + \chi_2(t,(t,\xi,\eta))
\end{array}
\right).}
$$

Thus $\bar{H}(t,\xi,\eta)$ is continuous on $\mathbb{R} \times \mathbb{R}^{n_1+n_2}$ and
$$
\left|
\bar{H}(t,\xi,\eta) -
{\left(
\begin{array}{c}
    \xi \\
    \eta
\end{array}
\right)}
\right|
\leq \bar{\sigma}.
$$

Moreover,
$$
\bar{H}(t,x(t,\tau,\xi,\eta),y(t,\tau,\xi,\eta)) =
{\left(
\begin{array}{c}
    x(t,\tau,\xi,\eta) + \chi_1(t,(t,x(t,\tau,\xi,\eta),y(t,\tau,\xi,\eta)))\\
    y(t,\tau,\xi,\eta) + \chi_2(t,(t,x(t,\tau,\xi,\eta),y(t,\tau,\xi,\eta)))
\end{array}
\right),}
$$
where
$
\chi(s,(t,x(t,\tau,\xi,\eta),y(t,\tau,\xi,\eta)))
$
is the unique bounded solution of system
$$
\frac{dz}{ds} =
W(s)z(s) +W_0(s)z(\gamma(s))+ \bar{h}(s,z(s),z(\gamma(s)),(t,x(t,\tau,\xi,\eta),y(t,\tau,\xi,\eta))).
$$

From
$$
x(s,(t,x(t,\tau,\xi,\eta),y(t,\tau,\xi,\eta))) = x(s,\tau,\xi,\eta),
$$
$$
y(s,(t,x(t,\tau,\xi,\eta),y(t,\tau,\xi,\eta))) = y(s,\tau,\xi,\eta),
$$
we have
$$\bar{h}(s,z(s),z(\gamma(s)),(t,x(t,\tau,\xi,\eta),y(t,\tau,\xi,\eta))) = \bar{h}(s,z(s),z(\gamma(s)),(\tau,\xi,\eta)).$$
Thus
$$\chi(s,(t,x(t,\tau,\xi,\eta),y(t,\tau,\xi,\eta))) = \chi(s,(\tau,\xi,\eta)), \quad \forall s \in \mathbb{R}.$$

Taking $s=t$, we get
$$\chi(t,(t,x(t,\tau,\xi,\eta),y(t,\tau,\xi,\eta))) = \chi(t,(\tau,\xi,\eta)).$$

Therefore,
$
\bar{H}(t,x(t,\tau,\xi,\eta),y(t,\tau,\xi,\eta)) =
{\left(
\begin{array}{c}
    x(t,\tau,\xi,\eta)+\chi_1(t,(\tau,\xi,\eta))\\
    y(t,\tau,\xi,\eta)+\chi_2(t,(\tau,\xi,\eta))
\end{array}
\right).}
$

We could check that $\bar{H}(t,x(t,\tau,\xi,\eta),y(t,\tau,\xi,\eta))$ is a solution of system \eqref{sysnlxy-pq} and \\
$|\bar{H}(t,x(t,\tau,\xi,\eta),y(t,\tau,\xi,\eta))-(x(t,\tau,\xi,\eta), y(t,\tau,\xi,\eta))^T|$ is bounded.  Therefore $\bar{H}(t,x,y)$ satisfies (i) and (ii).

Assume that $\bar{K}(t,x,y)={\left(
\begin{array}{c}
    \bar{K}_1(t,x,y) \\
    \bar{K}_2(t,x,y)
\end{array}
\right)}$ satisfies (\romannumeral1) and (\romannumeral2), too, where $\bar{K}_1(t,x,y)\in \mathbb{R}^{n_1}$ and $\bar{K}_2(t,x,y)\in \mathbb{R}^{n_2}$.
Since
$
{\left(
\begin{array}{c}
    x(t,\tau,\xi,\eta) \\
    y(t,\tau,\xi,\eta)
\end{array}
\right)}
$
is the solution of system \eqref{sysnlxy-phi}, $\bar{K}(t,x(t,\tau,\xi,\eta),y(t,\tau,\xi,\eta))$ is a solution of system \eqref{sysnlxy-pq}.

Denote
$
w(t) =
{\left(
\begin{array}{c}
    w_1(t) \\
    w_2(t)
\end{array}
\right)}
=
{\left(
\begin{array}{c}
    \bar{K}_1(t,x(t,\tau,\xi,\eta),y(t,\tau,\xi,\eta)) -x(t,\tau,\xi,\eta) \\
    \bar{K}_2(t,x(t,\tau,\xi,\eta),y(t,\tau,\xi,\eta))-y(t,\tau,\xi,\eta)
\end{array}
\right).}
$

From $w'(t)=W(t)w(t) +W_0(t)w(\gamma(t))+ \bar{h}(t,w(t),w(\gamma(t)),(\tau,\xi,\eta))$, we have that $w(t)$ is a bounded solution of system \eqref{Lemma7.1-temp1}. Therefore
$$w(t)=\chi(t,(\tau,\xi,\eta)).$$

Thus
$
\bar{K}(t,x(t,\tau,\xi,\eta),y(t,\tau,\xi,\eta)) =
{\left(
\begin{array}{c}
    x(t,\tau,\xi,\eta)+ \chi_1(t,(\tau,\xi,\eta)) \\
    y(t,\tau,\xi,\eta)+ \chi_2(t,(\tau,\xi,\eta))
\end{array}
\right).}
$

Taking
$
t = \tau$, we have
$$\bar{K}(\tau,\xi,\eta) =
{\left(
\begin{array}{c}
    \xi+ \chi_1(\tau, (\tau,\xi,\eta)) \\
    \eta+ \chi_2(\tau, (\tau,\xi,\eta))
\end{array}
\right)}
= \bar{H}(\tau,\xi,\eta).
$$

Thus $\bar{H}(t,x,y)$ is a unique function satisfying the conditions (\romannumeral1) and (\romannumeral2).  We complete the proof.   $\Box$

\th{Lemma 7.2} \emph{ System \eqref{sysnlxy-fg} is topologically conjugate to system \eqref{sysnlxy-phi}.  }

\proof
From Lemma 7.1, for any $t \in \mathbb{R}$,  $x, \tilde{x} \in \mathbb{R}^{n_1}$  and $y, \tilde{y} \in \mathbb{R}^{n_2}$, there exists a unique function $\tilde{H}(t,x,y)$ satisfies that

\begin{description}
\item[(i)]
There exists a constant $\sigma_1>0$ such that
$$
|\tilde{H}(t,x,y) - (x,y)^T|\leq \sigma_1.
$$

\item[(ii)] If
$
{\left(
\begin{array}{c}
    x(t) \\
    y(t)
\end{array}
\right)}
$
is a solution of system \eqref{sysnlxy-phi}, then $H(t,x(t),y(t))$ is a solution of system \eqref{sysnlxy-fg}.
\end{description}

Similarly, there exists a unique function $\tilde{L}(t,\tilde{x},\tilde{y})$ satisfies that

\begin{description}
\item[(i)]
There exists a constant $\sigma_2>0$ such that
$$
|\tilde{L}(t,\tilde{x},\tilde{y}) -(\tilde{x},\tilde{y})^T|\leq \sigma_2.
$$

\item[(ii)] If
$
{\left(
\begin{array}{c}
    \tilde{x}(t) \\
    \tilde{y}(t)
\end{array}
\right)}
$
is a solution of system \eqref{sysnlxy-fg}, then $\tilde{L}(t,\tilde{x}(t),\tilde{y}(t))$ is a solution of system \eqref{sysnlxy-phi}.
\end{description}

In what followings, we prove that
$
\tilde{L}(t,\tilde{H}(t,x,y)) =(x,y)^T$ and
$\tilde{H}(t,\tilde{L}(t,x,y)) =(x, y)^T$.

Denote $\tilde{J}(t,x,y)=\tilde{L}(t,\tilde{H}(t,x,y))$.

If
$
{\left(
\begin{array}{c}
    x(t) \\
    y(t)
\end{array}
\right)}
$
is a solution of system \eqref{sysnlxy-phi}, then $\tilde{H}(t,x(t),y(t))$ is a solution of system \eqref{sysnlxy-fg}. Thus $\tilde{L}(t,\tilde{H}(t,x(t),y(t)))$ is a solution of system \eqref{sysnlxy-phi}. By a simple calculation, we get
$$
|\tilde{J}(t,x,y) -(x, y)^T|\leq |\tilde{L}(t,\tilde{H}(t,x,y))-\tilde{H}(t,x,y)| +|\tilde{H}(t,x,y) -(x,y)^T|\leq \sigma_1+\sigma_2.
$$
Therefore $\tilde{J}(t,x,y)$ is the unique function satisfying the conditions  (\romannumeral1) and (\romannumeral2) in Lemma 7.1 which transforms the solution of system \eqref{sysnlxy-fg} to those of itself.

In particular, taking $p=\phi$ and $=\psi$ in system \eqref{sysnlxy-pq}, then  system \eqref{sysnlxy-pq} becomes system \eqref{sysnlxy-phi}. From system \eqref{sysnlxy-phi} to itself,  for any $t \in \mathbb{R}, x \in \mathbb{R}^{n_1}, y \in \mathbb{R}^{n_2}$,  the function
$\bar{H}(t,x,y)=
{\left(
\begin{array}{c}
    x \\
    y
\end{array}
\right)}
$
satisfies the conditions  (\romannumeral1) and (\romannumeral2)  in Lemma 7.1. Thus, for any $t \in \mathbb{R}$,  $x \in \mathbb{R}^{n_1}$ and $y \in \mathbb{R}^{n_2},$
$$
\tilde{J}(t,x,y)=\bar{H}(t,x,y)=
{\left(
\begin{array}{c}
    x \\
    y
\end{array}
\right)}.$$

That is
$$
\tilde{L}(t,\tilde{H}(t,x,y))=
{\left(
\begin{array}{c}
    x \\
    y
\end{array}
\right),}
\quad \forall t \in \mathbb{R}, x \in \mathbb{R}^{n_1}, y \in \mathbb{R}^{n_2}.
$$

Applying Lemma 7.1 to system \eqref{sysnlxy-pq} with $p=0$ and $q=0$, we could prove that
$$
\tilde{H}(t,\tilde{L}(t,\tilde{x},\tilde{y}))=
{\left(
\begin{array}{c}
    \tilde{x} \\
    \tilde{y}
\end{array}
\right),}
\quad \forall t \in \mathbb{R}, \tilde{x} \in \mathbb{R}^{n_1}, \tilde{y}\in \mathbb{R}^{n_2}.
$$

Therefore, for a fixed $t$, $\tilde{H}^{-1}(t,\cdot,\cdot)=\tilde{L}(t,\cdot,\cdot)$.

According to Lemma 5.4 and Remark 5.3, solutions of systems \eqref{sysnlxy-pq} and  \eqref{sysnlxy-phi}  are continuous with respect to initial values.

Since both $\tilde{H}(t,\cdot)$ and $\tilde{L}(t,\cdot)$ are continuous, $\tilde{H}(t,\cdot)$ and $\tilde{L}(t,\cdot)$ are homeomorphisms of $\mathbb{R}^n$. Thus System \eqref{sysnlxy-fg} is topologically conjugate to system \eqref{sysnlxy-phi}.  The proof is complete.   $\Box$

\section{The proof of Theorem 2}
From Lemmas 6.12 and  7.2, we have that $H(t,\cdot)\circ\tilde{H}(t,\cdot)$ and $\tilde{L}(t,\cdot)\circ L(t,\cdot)$ are homeomorphisms of $\mathbb{R}^n$ and $\big{(}H(t,\cdot)\circ\tilde{H}(t,\cdot)\big{)}^{-1}=\tilde{L}(t,\cdot)\circ L(t,\cdot)$. Moreover,  $H(t,\cdot)\circ\tilde{H}(t,\cdot)$ sends the solutions of system \eqref{syslxy} onto those of system \eqref{sysnlxy-phi} and $\tilde{L}(t,\cdot)\circ L(t,\cdot)$ sends the solutions of system \eqref{sysnlxy-phi} onto those of system \eqref{syslxy}. It is easy to see that $|H(t,\cdot)\circ\tilde{H}(t,(x,y)^T)-(x,y)^T|$ and $|\tilde{L}(t,\cdot)\circ L(t,(x,y)^T)-(x,y)^T|$ are bounded. Therefore system \eqref{sysnlxy-phi} and system \eqref{syslxy} are topologically conjugated. $\Box$

\section{Conflict of Interests}

The authors declare that there is no conflict of interests
regarding the publication of this article.


\end{document}